\documentclass[12pt]{amsart}

\usepackage{amssymb}
\usepackage{latexsym}
\usepackage{amsmath}

\begin{document}

\def\l{{\mathcal L}}
\def\p{{\mathcal P}}
\def\c{{\mathbb C}}
\def\r{{\mathbb R}}
\def\h{{\mathbb H}}

\newcommand\wt{\widetilde}
\newcommand\wh{\widehat}
\newcommand\wo{\overline} 

\newcommand{\e}[2]{\begin{pmatrix}{#2}\cr{#1}\end{pmatrix}}
\newcommand{\ee}[4]{\begin{pmatrix}{#2}&{#4}\cr{#1}&{#3}\end{pmatrix}}
\newcommand{\eee}[6]{\begin{pmatrix}{#2}&{#4}&{#6}\cr{#1}&{#3}&{#5}\end{pmatrix}}

\def\wtm{\,\widetilde{\ }\,}
\def\whm{\,\widehat{\ }\,}

\newtheorem{theo}{Theorem}[section]
\newtheorem{prop}{Proposition}[section]
\newtheorem{lemm}{Lemma}[section]
\newtheorem{defi}{Definition}[section]
\newtheorem{coro}{Corollary}[section]
\newtheorem{exam}{Example}[section]
\newtheorem{rema}{Remark}[section]

\title{The planar algebra of a coaction}
\author{Teodor Banica}
\address{Institut de Mathematiques de Jussieu, 175 rue du Chevaleret, 75013 Paris}
\email{banica@math.jussieu.fr}

\maketitle

{\bf Abstract.} We study actions of ``compact quantum groups'' on ``finite quantum spaces''. According to Woronowicz and to general $\c^*$-algebra philosophy these correspond to certain coactions $v:A\to A\otimes H$. Here $A$ is a finite dimensional $\c^*$-algebra, and $H$ is a certain special type of Hopf $*$-algebra. If $v$ preserves a positive linear form $\varphi :A\to\c$, a version of Jones' ``basic construction'' applies. This produces a certain $\c^*$-algebra structure on $A^{\otimes n}$, plus a coaction $v_n :A^{\otimes n}\to A^{\otimes n}\otimes H$, for every $n$. The elements $x$ satisfying $v_n(x)=x\otimes 1$ are called fixed points of $v_n$. They form a $\c^*$-algebra $Q_n(v)$. We prove that under suitable assumptions on $v$ the graded union of the algebras $Q_n(v)$ is a spherical $\c^*$-planar algebra.

\section*{Introduction}

A major discovery of the eighties, due to Jones \cite{j0}, is that an inclusion of von Neumann algebras with trivial centers, also called subfactor, produces a representation of the algebra of Temperley and Lieb \cite{tl}. This can be used for getting several unexpected results about von Neumann algebras. For instance that subfactors of index $< 4$ are classified by ADE diagrams, and that their indices must be of the form $4\cos^2(\pi /n)$.

The classification program for subfactors, with many people involved over the last 20 years, already reached a few final conclusions. Among them is an axiomatization of a large class of ``quantum algebras'', having positivity properties. A first set of axioms, of algebraic nature, was found by Popa in \cite{po}. A set of topological axioms, leading to the notion of planar algebra, was found by Jones in \cite{j1}.

The colored planar operad $\p$ consists of certain planar diagrams called tangles. Each tangle has several input discs and an output disc, connected by non-crossing strings. The operad law is given by gluing of tangles. A planar algebra is by definition an algebra over $\p$. That is, we have a graded vector space $Q=Q_0^\pm,\, Q_1,\, Q_2,\, Q_3,\ldots$, and for every tangle $T$ we can put elements of $Q$ in the input discs of $T$ and we get an element of $Q$ on the output disc of $T$.

When the scalars are complex numbers and certain positivity and spherical invariance properties hold, $Q$ is said to be spherical $\c^*$-planar algebra. Results of Jones \cite{j1} and Popa \cite{po} show that every subfactor produces a such a planar algebra, and vice versa. In the ``amenable'' case the correspondence is one-to-one, by a result of Popa \cite{p0}.

\medskip

It is natural to ask about how these fundamental techniques from subfactors work for compact quantum groups. According to Woronowicz \cite{w3} such a quantum group is described by a certain special type of Hopf $\c^*$-algebra. So, let $H$ be such a Hopf $\c^*$-algebra, let $A$ be a finite-dimensional $\c^*$-algebra and let $v :A\to A\otimes H$ be a coaction. It is convenient to assume that $v$ is co-faithful, in the sense that its coefficients generate $H$ as a $\c^*$-algebra.

If $v$ leaves invariant a linear form the basic construction produces coactions $v_n :A^{\otimes n}\to A^{\otimes n}\otimes H$ for every $n$. Here the tensor powers $A^{\otimes n}$ are given the $\c^*$-algebra structure coming from the basic construction. Consider the algebras $Q_n(v)$ of fixed points under the coactions $v_n$. That is, of elements satisfying $v_n(x)=x\otimes 1$. These form an increasing sequence of finite dimensional $\c^*$-algebras, and their union is a graded $*$-algebra, denoted $Q(v)$.

\medskip

In most cases of interest $Q(v)$ is known to be a spherical $\c^*$-planar algebra. Moreover, $Q(v)$ encodes important information about $(H,v)$, and several algebraic or analytic properties of $(H,v)$ can be translated in terms of $Q(v)$. In fact, it is expected that a reconstruction map of type $Q(v)\to (H,v)$ exists, as a modification of Woronowicz's Tannakian duality \cite{w2}.

When $H$ is finite dimensional and $v:H\to H\otimes H$ is its comultiplication, this follows from Ocneanu's depth 2 duality, see David \cite{da}, Longo \cite{lo} and Szymanski \cite{sz}. A direct proof is obtained by Kodiyalam, Landau and Sunder in \cite{kls}. For more results on the depth 2 case see Das \cite{das} and Das and Kodiyalam \cite{dak}.

More generally, one can consider the case when $H$ is a Kac type, meaning that the square of its antipode $S^2$ is the identity. Several explicit results, due to Landau \cite{la}, Landau and Sunder \cite{ls}, and Bhattacharyya and Landau \cite{bl} are available here. In the general $S^2=id$ case a subfactor is constructed in \cite{kac}, and its standard invariant is computed by using a method of Wassermann from \cite{was}. By combining this with a result of Jones in \cite{j1}, it follows that $Q(v)$ is a planar algebra.

In the $S^2\neq id$ case things are less explicit. When $H=C(G)_q$ corresponds to a $q$-deformation with $q>0$ of a compact Lie group and $v$ comes from a projective representation of $G$ this follows from work of Sawin \cite{sa}. More generally, when $A=M_n(\c )$ and $v$ is adjoint to a corepresentation of $H$, this follows from a many-to-one Tannakian correspondence, established in \cite{subf}.

\medskip

The problem with most of the above results is that the planar algebra structure of $Q(v)$ is not quite explicit, because it comes from a subfactor or a standard $\lambda$-lattice, via the fundamental results of Jones \cite{j1} and Popa \cite{po}. The other obvious problem is that all these results certainly cover the most interesting cases, but some cases are still left. And finally, a third problem is with the reconstruction map, not available in most cases.

One may wonder about a very general correspondence of the form $(H,A,v)\leftrightarrow Q(v)$, between triples $(H,A,v)$ satisfying certain assumptions and certain planar algebras. Something like ``triples satisfying a Perron-Frobenius type condition are in one-to-one correspondence with twisted $\c^*$-subalgebras of depth 1 planar algebras''. Moreover, for this result to be ready to use, one would like to have a direct construction of the correspondence, somehow in the spirit of Ocneanu's depth 2 duality and of Woronowicz's Tannakian duality.

So far, the only fully satisfactory result in this sense seems to be the one in the depth 2 case, where the enlightening paper of Kodiyalam, Landau and Sunder \cite{kls} is available.

The aim of the present work is to construct a general map of type $(H,A,v)\to Q(v)$.

In \S 1 and \S 2 we apply the basic construction, and we study the equivariance properties of various annular tangles. In the planar algebra setting it is convenient to use bases and indices and wo do it right from the beginning. This requires a normalisation of the coefficients of $v$. We choose the one which makes the spin factor behave uniformly at even and odd levels.

In \S 3 and \S 4 we prove that under suitable assumptions $Q(v)$ is a spherical $\c^*$-planar algebra. When the square of the antipode $S^2$ is the identity this is a subalgebra of the depth 1 planar algebra $P(A)$ constructed by Jones in \cite{j2}, by using a certain explicit statistical mechanical sum. In the general case the inclusion $Q(v)\subset P(A)$ appears to be ``twisted'', and the partition function of $Q(v)$ comes here from a standard $\lambda$-lattice in the sense of Popa \cite{po}, by using the ``bubbling''construction of Jones \cite{j1}.

\medskip

As a conclusion, in the $S^2=id$ case the map $(H,A,v)\to Q(v)$ is constructed quite explicitely, and what is left is to do the converse construction. In the $S^2\neq id$ case what we do is rather to compute the domain of $(H,A,v)\to Q(v)$, by a method which is to be improved.

The first version of this paper was written in 2002. This version is the third one, written in 2004, with new notations and many comments added, but basically containing the same material. So far, we have found no improvement in the $S^2\neq id$ case. 

In the recent paper \cite{gr} we obtain the duality for coactions on $A=\c^n$. Here the condition $S^2=id$ is automatic. This duality restricts to a correspondence between Hopf $\c^*$-algebras associated to colored graphs with $n$ vertices and planar subalgebras of the spin planar algebra $P(\c^n)$, generated by a self-adjoint 2-box. This latter correspondence makes a link between Hopf $\c^*$-algebras and the classification program initiated by Bisch and Jones in \cite{bj1} and \cite{bj2}, and can be used for explicit (numeric) computations of Poincar\'e series of such Hopf $\c^*$-algebras.

Some other possible applications of such dualities are discussed in \S 5 in \cite{gr}.

\section{Formalism}

We are interested in coactions of the form $v:A\to A\otimes H$, where $A$ is a finite dimensional $\c^*$-algebra, and $H$ is a Hopf $*$-algebra corresponding to a ``compact quantum group''. For instance we will need at some point the existence of a Haar functional $h:H\to\c$.

Moreover, we want to deal with the general case, where the square of the antipode is not necessarely the identity, $S^2\neq id$. This is known to correspond to the case where $h$ has not the trace property, $h(ab)\neq h(ba)$. This is why we will also need at some point a precise description of the modular automorphism of $h$. That is, of the automorphism $\sigma :H\to H$ satisfying $h(ab)=h(b\sigma (a))$ for any $a,b$.

So, what we need is the notion of Hopf $*$-algebra associated to a compact matrix pseudogroup, in the sense of Woronowicz \cite{w1}. Indeed, for such a Hopf $*$-algebra, both the existence of $h$ and the explicit computation of $\sigma$ appear as theorems. Conversely, given a Hopf $*$-algebra with Haar integral etc. coacting via $v:A\to A\otimes H$ on a finite dimensional $\c^*$-algebra, the coaction $v$ can be regarded as a corepresentation of $H$ on the linear space $A$. Nothing is lost when replacing $H$ with the $*$-algebra generated by the coefficients of $v$, and $v$ becomes in this way a fundamental corepresentation of $H$ in the sense of \cite{w1}.

However, it is more convenient to start with the more enligthening axioms in Woronowicz's paper \cite{w3}. A good reference here is the paper \cite{mvd} by Maes and Van Daele, containg a short exposition of the subject, with several simplifications, and available at arxiv.org.

The terminology in the definition below is probably quite reasonable, but not standard.

\begin{defi}
A Hopf $\c^*$-algebra with unit is a pair $\h =(\h,\Delta )$ consisting of a $\c^*$-algebra with unit $\h$ and a $\c^*$-morphism $\Delta :\h\to \h\otimes \h$, subject to the following conditions.

(i) Coassociativity condition $(id\otimes\Delta )\Delta =(\Delta \otimes id )\Delta$.

(ii) Cocancellation law: the sets $\Delta (\h)(1\otimes \h)$ and $\Delta (\h)(\h\otimes 1)$ are dense in $\h\otimes \h$.
\end{defi}

The basic example is $\h=\c (G)$, the algebra of continuous functions on a compact group $G$, with $\Delta (\varphi): (g,h)\to\varphi (gh)$. Here coassociativity of $\Delta$ follows from associativity of the multiplication $\cdot$ of $G$, and cocancellation in $(\h,\Delta )$ follows from cancellation in $(G,\cdot )$.

Conversely, assume that $(\h,\Delta )$ is as in definition 1.1, and that $\h$ is commutative. The Gelfand transform gives an isomorphism $\h\simeq \c (G)$, where $G$ is the spectrum of $H$. Now the coassociative map $\Delta$ gives rise to an associative map $\cdot :G\times G\to G$. In other words, we have here a compact semigroup $(G,\cdot )$, which by (ii) follows to have cancellation. It is then well-known that $G$ must be a compact group.

As a conclusion, the construction $(G,\cdot )\to (\c (G),\Delta )$ is a contravariant equivalence of categories between compact groups and commutative Hopf $\c^*$-algebras with unit. So, a pair $(\h,\Delta )$ as in definition 1.1 can be thought of as corresponding to a ``compact quantum group''.

Among main results of Woronowicz in \cite{w3} is the construction of a dense subalgebra $H\subset\h$, consisting of ``representative functions'' on the compact quantum group. This has a counit $\varepsilon :H\to\c$ and an antipode $S:H\to H$, which satisfy the usual Hopf algebra identities. It is convenient to dentode by $m :H\otimes H\to H$ and $u:\c\to H$ its multiplication and unit maps.

\begin{defi}
In this paper $H=(H,m,u,\Delta ,\varepsilon, S,*)$ will denote the Hopf $*$-algebra of ``representative functions on a compact quantum group'', in the sense that $H$ is the canonical dense subalgebra associated by Woronowicz to a Hopf $\c^*$-algebra with unit $\h$.
\end{defi}

As already explained, $H$ will come in fact together with a fundamental corepresentation. So, we will feel free to refer to results of Woronowicz from the fundamental paper \cite{w1}.

The other piece of data we need is a pair $(A,\varphi )$ consisting of a finite dimensional $\c^*$-algebra $A$ and a positive linear form $\varphi :A\to\c$. It is well-known that $A$ must be isomorphic to a direct sum of matrix algebras, and $\varphi$ must be of the form $a\to tr(qa)$, with $q\in A$ positive.

A basic example here is the algebra $A=\c (X)$ of functions on a finite set $X$, with linear form $\varphi (f)=\sum f(x)\mu (x)$, where $\mu$ is a positive measure on $X$. It is probably tempting to think of a general pair $(A,\varphi )$ as corresponding to a ``measured finite quantum space''. But the other main example is $A=M_n (\c)$ with $\varphi (a)=tr(qa)$, where $q$ is some positive matrix, and here this interpretation doesn't quite help.

\begin{defi}
Let $(A,\varphi )$ be a finite dimensional $\c^*$-algebra together with a positive linear form. A coaction of $H$ on $(A,\varphi )$ is a morphism of $*$-algebras $v :A\to A\otimes H$, subject to the following conditions.

(i) Coassociativity condition $(v\otimes id)v =(id\otimes\Delta )v$.

(ii) Counitality condition $(id\otimes\varepsilon )v =id$.

(iii) Copreservation of $\varphi$ condition $(\varphi\otimes id)v =\varphi (.)1$.
\end{defi}

If (iii) is not satisfied we just say that $v$ is a coaction of $H$ on $A$.

The purpose of this section is to reformulate these axioms, for further use in establishing results about propagation of $v$ in the Jones tower for $\c\subset A$. The precise structure of the Jones tower for $\c\subset A$ is that of a spherical $\c^*$-planar algebra. The following are known.

(1) Bases and indices are needed so far in understanding this planar algebra structure, meaning that an approach with ``global formulae'' is not available yet. In fact, a planar algebra is quite an abstract notion, and the action of tangles on tensors is best understood by keeping in mind rules like ``indices are allowed to travel on strings'' or ``two different indices make the whole thing vanish when they meet'' etc. This is why indices are necessary.

(2) Some quite unobvious choices of bases, normalisations, notations etc. are needed as well. See e.g. the comments of Jones in \cite{j1} and \cite{j2}. The idea here is that the planar meaning of various ``deformation'' parameters is very unclear. The ``spin vector'' used by Jones in \cite{j2}, which already requires a tricky normalisation, turns to have a quite clear planar interpretation, in terms of ``horizontal'' structure. In this paper the set of ``parameters'' will be even bigger. This will require several careful normalisations, and the problem of finding a reasonable planar interpretation of these parameters will be eventually left open in the general case. 

Now (1) tells us to look for a reformulation of definition 1.3, in terms of coefficients of $v$, with respect to some basis of $A$. This is an a priori quite standard task: coassociativity corresponds to the well-known condition $\Delta (v_{ij})=\sum v_{ik}\otimes v_{kj}$ and so on. However, because of (2), we have be extremely careful in the choice of the basis and coefficients.

We will use a normalisation which may seem a bit strange, but which does work, in the sense that formulae in the Jones tower will look quite similar at even and odd levels. Of course, this choice of simplifying things in higher formulae to come might cause the very first formulae -- in statements and proofs -- to look more complicated than needed. This will be indeed the case.

\begin{defi}
Let $(A,\varphi )$ be as above. Choose a system of matrix units $X\subset A$ making $\varphi$ diagonal, with the following multiplication convention.
$$\e{i}{j}\e{k}{l}=\delta_{jk}\e{i}{l}$$

We denote by $q_i$ the fourth roots of the weights of $\varphi$, chosen positive.
$$\varphi\e{i}{j}=\delta_{ij}q_i^4\hskip 2cm q_i> 0$$

Any linear map $v:A\to A\otimes H$ will be written in the following form.
$$v\e{i}{j}=\sum\e{k}{l}\otimes 
q_k^{-1}q_iq_jq_l^{-1}V\ee{k}{l}{i}{j}$$

This is, to any linear map $v$ we associate in this way a matrix $V$, and vice versa.
\end{defi}

It is convenient to define the coefficients $V(^l_k{\ }^j_i)$ for all indices $i,j,k,l$, by saying that they are equal to zero if $(^j_i)$ or $(^l_k)$ don't exist. In fact, best here would be to use the groupoid structure os $X$, but since we don't have results for more general groupoids, we don't do it.

As for the sum sign in definition 1.4, this is by definition over all elements $(^l_k)\in X$. More generally, in any formula of type $A=\sum B$ or $\sum B=A$ with $A,B\in H$ the sum will be over all indices which appear in $B$ and don't appear in $A$.

The normalisation in definition 1.4 is the one which will appear to work well in the Jones tower. For, we must first do the above-mentioned reformulation of definition 1.3.

\begin{prop}
A linear map $v :A\to A\otimes H$ is a coaction of $H$ on $A$ if and only if $V$ satisfies the following conditions.
$$\varepsilon{V}\ee{k}{l}{i}{j} =\delta_{ki}\delta_{lj}1$$
$$\Delta{V}\ee{k}{l}{i}{j} =\sum {V}\ee{k}{l}{g}{h}\otimes{V}\ee{g}{h}{i}{j}$$
$$V\ee{k}{l}{i}{j}^*=V\ee{l}{k}{j}{i}$$ 
$$\sum q_i^2V\ee{k}{l}{i}{i}=\delta_{kl}q_k^2$$ 
$$\sum q_s^{-2}{V}\ee{k}{s}{g}{h}{V}\ee{s}{l}{i}{j}=\delta_{hi}q_i^{-2}{V}\ee{k}{l}{g}{j}$$

This sequence of five conditions will be denoted ($\varepsilon$), ($\Delta$), ($*$), ($u^\circ$), (${\, }^\circ m$).
\end{prop}

\begin{proof}
This is well-known, modulo our normalisations for $V$, so the only thing to check is that all $q$ values in the statement are the good ones. It is possible to prove this either by using global formulae, or with a direct matrix computation. We prefer to present this latter approach, as a warm-up for more involved computations to come, where bases and matrix computations seem to be unavoidable, cf. above considerations (1) and (2).

By using the defining formula of $v$ we get
$$(id\otimes\varepsilon )v\e{i}{j} =\sum \e{k}{l}\otimes 
q_k^{-1}q_iq_jq_l^{-1}\varepsilon V\ee{k}{l}{i}{j}$$
so the condition $(id\otimes\varepsilon )v =id$ holds if and only if $V$ satisfies
$$\varepsilon V\ee{k}{l}{i}{j}=
q_kq_i^{-1}q_j^{-1}q_l\delta_{ki}\delta_{lj}$$
for any $i,j,k,l$, i.e. if and only if $V$ satisfies ($\varepsilon$). We have
\begin{eqnarray*}
(v\otimes id)v \e{i}{j}
&=&\sum v\e{k}{l}\otimes q_k^{-1}q_iq_jq_l^{-1}V\ee{k}{l}{i}{j}\\ 
&=&\sum\e{g}{h}\otimes q_g^{-1}q_kq_lq_h^{-1}V\ee{g}{h}{k}{l}\otimes
q_k^{-1}q_iq_jq_l^{-1}V\ee{k}{l}{i}{j}\\
&=&\sum\e{g}{h}\otimes q_g^{-1}q_h^{-1}q_iq_jV\ee{g}{h}{k}{l}\otimes V\ee{k}{l}{i}{j}\\
&=&\sum\e{k}{l}\otimes q_k^{-1}q_l^{-1}q_iq_jV\ee{k}{l}{g}{h}\otimes V\ee{g}{h}{i}{j}
\end{eqnarray*}
so $(v\otimes id)v =(id\otimes\Delta )v$ is equivalent to ($\Delta$). We have
$$v (1)=\sum v\e{i}{i}=\sum \e{k}{l}\otimes 
q_k^{-1}q_i^2q_l^{-1}V\ee{k}{l}{i}{i}$$
so $v (1)=1\otimes 1$ is equivalent to ($u^\circ$). Also, from the formulae
\begin{eqnarray*}
v\e{g}{h}v\e{i}{j}
&=&\sum\e{k}{s}\e{S}{l}\otimes q_k^{-1}q_gq_hq_s^{-1}q_S^{-1}q_iq_jq_l^{-1}
V\ee{k}{s}{g}{h}V\ee{S}{l}{i}{j}\\
&=&\sum\e{k}{l}\otimes (q_i^2q_s^{-2})(q_hq_i^{-1})
q_k^{-1}q_gq_jq_l^{-1}V\ee{k}{s}{g}{h}V\ee{s}{l}{i}{j}\\
v\left( \e{g}{h}\e{i}{j}\right)
&=&\sum\e{k}{l}\otimes \delta_{hi}q_k^{-1}q_gq_jq_l^{-1}V\ee{k}{l}{g}{j}
\end{eqnarray*}
we get that $v$ is multiplicative if and only if (${\, }^\circ m$) holds.We have
$$v\e{j}{i}=\sum\e{l}{k}\otimes 
q_l^{-1}q_jq_iq_k^{-1}V\ee{l}{k}{j}{i}$$
so $v$ is involutive if and only if ($*$) holds.
\end{proof}

\begin{prop}
Assume that $v$ is a coaction of $H$ on $A$. Then the following three conditions ($S$), (${\, }^\circ u$) and ($m^\circ$)
$$S{V}\ee{k}{l}{i}{j} =q_k^2q_i^{-2}q_j^2q_l^{-2}{V}\ee{j}{i}{l}{k}$$
$$\sum q_i^2V\ee{i}{i}{k}{l}=\delta_{kl}q_k^2$$
$$\sum q_s^{-2}V\ee{k}{h}{g}{s}{V}\ee{i}{l}{s}{j}=\delta_{hi}q_i^{-2}{V}\ee{k}{l}{g}{j}$$
are equivalent, and are satisfied if and only if $v$ preserves $\varphi$.
\end{prop}

\begin{proof}
We keep proving things by performing matrix computations. We have
$$(\varphi\otimes id)v \e{k}{l}
=\sum \varphi\e{i}{j}q_i^{-1}q_kq_lq_j^{-1} V\ee{i}{j}{k}{l}\\
=q_kq_l\sum q_i^2 V\ee{i}{i}{k}{l}$$
so the condition $(\varphi\otimes id)v =\varphi (.)1$ holds if and only if $V$ satisfies (${\, }^\circ u$). It remains to prove that if ($\varepsilon$), ($\Delta$), ($u^\circ$), (${\, }^\circ m$) are satisfied, then ($S$), ($m^\circ$) and (${\, }^\circ u$) are equivalent. First, by applying $S$ to ($u^\circ$) we get that  ($S$) implies (${\, }^\circ u$).
$$\delta_{kl}q_k^2=\sum q_i^2SV\ee{k}{l}{i}{i}=
\sum q_i^2q_k^2q_l^{-2}V\ee{i}{i}{k}{l}$$

Assume that (${\, }^\circ u$) holds. By combining it with (${\, }^\circ m$) we get
$$\sum q_i^2q_j^{-2}V\ee{i}{j}{k}{l}{V}\ee{j}{i}{g}{h}=
\sum q_i^2\delta_{lg}q_g^{-2}V\ee{i}{i}{k}{h}=\delta_{lg}q_g^{-2}\delta_{kh}q_k^2$$
and this can be rewritten in the following form.
$$\sum q_i^2q_k^{-2}q_l^2q_j^{-2}V\ee{i}{j}{k}{l}{V}\ee{j}{i}{g}{h}=\delta_{lg}\delta_{kh}1$$

If $e_{xy}$ with $x,y\in X$ is the system of matrix units in $\l (A)$ we get
$$\left( \sum e_{(^k_l)(^i_j)}\otimes q_i^2q_k^{-2}q_l^2q_j^{-2}V\ee{i}{j}{k}{l}
\right) \left( \sum e_{(^i_j)(^h_g)}\otimes {V}\ee{j}{i}{g}{h}\right)=1\otimes 1$$

On the other hand, conditions ($\varepsilon$) and ($\Delta$) say that $V\in\l (A)\otimes H$ is a corepresentation, i.e. that it satisfies 
$$(id\otimes\Delta )V=V_{12}V_{13}\hskip 2cm (id\otimes\varepsilon )V=1$$
so by considering $(id\otimes E)V$ with $E$ given by the Hopf algebra axiom
$$E=m(S\otimes id)\Delta =m(id\otimes S)\Delta =\varepsilon (.)1$$
we get that $(id\otimes S)V$ is an inverse for $V$. Thus the above formula gives $(id\otimes S)V$, and by identifying coefficients we get ($S$). Thus (${\, }^\circ u$) implies ($S$). It remains to prove that ($S$) is equivalent to ($m^\circ$). Assume that ($S$) is satisfied. Then 
\begin{eqnarray*}
&&S\left( \sum q_s^{-2}{V}\ee{k}{s}{g}{h}{V}\ee{s}{l}{i}{j}\right)\\
&=&\sum q_s^{-2}SV\ee{s}{l}{i}{j}SV\ee{k}{s}{g}{h}\\
&=&\sum q_s^{-2}q_s^2q_i^{-2}q_j^2q_l^{-2}q_k^2q_g^{-2}q_h^2q_s^{-2}
V\ee{j}{i}{l}{s}V\ee{h}{g}{s}{k}
\end{eqnarray*}

On the other hand by applying $S$ to the right term of (${\, }^\circ m$) we get
$$S\left( \delta_{hi}q_i^{-2}{V}\ee{k}{l}{g}{j}\right)
=\delta_{hi}q_i^{-2}q_k^2q_g^{-2}q_j^2q_l^{-2}V\ee{j}{g}{l}{k}$$

By using (${\, }^\circ m$) we get after cancelling $q$'s that
$$\delta_{hi}V\ee{j}{g}{l}{k}=
\sum q_h^2q_s^{-2}V\ee{j}{i}{l}{s}V\ee{h}{g}{s}{k}$$
and this is ($m^\circ$). Finally, the proof of ($m^\circ$) implies ($S$) is similar to the proof of (${\, }^\circ u$) implies ($S$). Indeed, by combining ($m^\circ$) and ($u^\circ$) we get
$$\sum q_i^2q_k^{-2}q_l^2q_j^{-2}V\ee{i}{j}{k}{l}{V}\ee{j}{i}{g}{h}=\delta_{lg}\delta_{kh}1$$
and this gives a right inverse for $V$, hence the formula ($S$) for the antipode.
\end{proof}

\section{Basic construction, equivariance results, and the untwisted case}

For any $n$ the set $X^n$ is a basis of the linear space $A^{\otimes n}$. With ``loop'' notations
$$\e{i_1}{j_1}\otimes\e{i_2}{j_2}\otimes\ldots\otimes\e{i_{2s}}{j_{2s}}=
\begin{pmatrix}j_{2s} & i_{2s} &\ldots & i_{s+1}\cr i_1 & j_1 &\ldots & j_s\end{pmatrix}$$
$$\e{i_1}{j_1}\otimes\e{i_2}{j_2}\otimes\ldots\otimes\e{i_{2s-1}}{j_{2s-1}}=
\begin{pmatrix}j_{2s-1} & i_{2s-1} &\ldots & j_s\cr i_1 & j_1 &\ldots & i_s\end{pmatrix}$$
for this basis, depending on the parity of $n$, the linear extension of
$$\begin{pmatrix}j_1 & \ldots & j_n\cr i_1 & \ldots & i_n\end{pmatrix}
\begin{pmatrix}l_1 & \ldots & l_n\cr k_1 & \ldots & k_n\end{pmatrix}=
\delta_{j_1k_1}\ldots \delta_{j_nk_n}
\begin{pmatrix}l_1 & \ldots & l_n\cr i_1 & \ldots & i_n\end{pmatrix}$$
is an associative multiplication on $A^{\otimes n}$. Together with the antilinear extension of
$$\begin{pmatrix}j_1 & \ldots & j_n\cr i_1 & \ldots & i_n\end{pmatrix}^*
=\begin{pmatrix}i_1 & \ldots & i_n\cr j_1 & \ldots & j_n\end{pmatrix}$$
this gives a finite dimensional $\c^*$-algebra structure on $A^{\otimes n}$. Note that $A^{\otimes 1}=A$ and that $A^{\otimes 2}$ is a matrix algebra. In fact the algebras $A^{\otimes n}$ are obtained from $A$ by performing the basic construction to the inclusion $\c\subset A$. See the book of Goodman, de la Harpe and Jones \cite{ghj}.

We use for $A^{\otimes n}$ the same conventions for sums etc. as those for $A$.

Conditions ($\varepsilon$) and ($\Delta$) show that the matrix
$$u=\sum e_{(^l_k)(^j_i)}\otimes V\ee{k}{l}{i}{j}\in\l (A)\otimes H$$
is a corepresentation, so we can consider its tensor powers.
$$u^{\otimes n}=u_{1,n+1}u_{2,n+1}\ldots u_{n,n+1}\in \l (A^{\otimes n})\otimes H$$

Let $V_n$ be the matrix of coefficients of $u^{\otimes n}$, defined by
$$u^{\otimes n}=\sum 
e_{(^{l_1\ldots l_n}_{k_1\ldots k_n})(^{j_1\ldots j_n}_{i_1\ldots i_n})}
\otimes V_n\begin{pmatrix}l_1&\ldots &l_n &j_1 &\ldots &j_n\cr 
k_1&\ldots &k_n &i_1 &\ldots &i_n\end{pmatrix}$$

Define a linear form $\wt{\varphi}_n$ by
$$\wt{\varphi}_n\begin{pmatrix}j_1 & \ldots & j_n\cr i_1 & \ldots & i_n\end{pmatrix}=
\delta_{(i_1\ldots i_n)(j_1\ldots j_n)}q_{(i_1\ldots i_n)}^4$$
where the weights are given by the function
$$q_{(i_1\ldots i_n)}=q_{i_1}q_{i_2}^{-1}q_{i_3}\ldots q_{i_n}^{\mp 1}$$
where $\pm 1=(-1)^n$. Note that $q_{(i)}=q_i$, so $\wt{\varphi}_1=\varphi$ on $A^{\otimes 1}=A$.

\begin{prop}
If $v :A\to A\otimes H$ is a coaction of $H$ on $A$ which preserves $\varphi$ then the linear map $v_n:A^{\otimes n}\to A^{\otimes n}\otimes H$ given by
$$v_n\begin{pmatrix}j_1 &\ldots &j_n\cr i_1 &\ldots &i_n\end{pmatrix}
=\sum \begin{pmatrix}l_1 &\ldots &l_n\cr k_1 &\ldots &k_n\end{pmatrix}
\otimes q_{(k_1\ldots k_n)}^{-1}q_{(i_1\ldots i_n)}
q_{(j_1\ldots j_n)}q_{(l_1\ldots l_n)}^{-1}$$
$$\cdot\,\, V_n\begin{pmatrix}l_1&\ldots &l_n &j_1 &\ldots &j_n\cr 
k_1&\ldots &k_n &i_1 &\ldots &i_n\end{pmatrix}$$
is a coaction of $H$ on $A^{\otimes n}$ which preserves $\wt{\varphi}_n$.
\end{prop}

\begin{proof}
The tensor powers of $u$ are given by
\begin{eqnarray*}
u^{\otimes n}
&=& \sum e_{(^{l_1}_{k_1})(^{j_1}_{i_1})}\otimes\ldots
\otimes e_{(^{l_n}_{k_n})(^{j_n}_{i_n})}\otimes V\ee{k_1}{l_1}{i_1}{j_1}\ldots 
V\ee{k_n}{l_n}{i_n}{j_n}\\
&=& \sum e_{(^{l_1}_{k_1})\otimes\ldots\otimes (^{l_n}_{k_n}),\, 
(^{j_1}_{i_1})\otimes\ldots\otimes (^{j_n}_{i_n})}\otimes 
V\ee{k_1}{l_1}{i_1}{j_1}\ldots V\ee{k_n}{l_n}{i_n}{j_n}\\
&=& \sum e_{(^{k_2}_{k_1})\otimes\ldots\otimes (^{l_1}_{l_2}),\, 
(^{i_2}_{i_1})\otimes\ldots\otimes (^{j_1}_{j_2})}\otimes 
V\ee{k_1}{k_2}{i_1}{i_2}\ldots V\ee{l_2}{l_1}{j_2}{j_1}\\
&=& \sum e_{(^{l_1\ldots l_n}_{k_1\ldots k_n})(^{j_1\ldots j_n}_{i_1\ldots i_n})}
\otimes V\ee{k_1}{k_2}{i_1}{i_2}\ldots V\ee{l_2}{l_1}{j_2}{j_1}
\end{eqnarray*}
and this gives the formula of $V_n$. More precisely, we have $V_1=V$ and
$$V_2\begin{pmatrix}l_1 & l_2 & j_1 & j_2 \cr k_1 & k_2 &i_1 & i_2\end{pmatrix}=
V\ee{k_1}{k_2}{i_1}{i_2}V\ee{l_2}{l_1}{j_2}{j_1}$$
$$V_3\begin{pmatrix}l_1 & l_2 &l_3& j_1 & j_2&j_3 \cr k_1 & k_2&k_3 &i_1 & i_2&i_3\end{pmatrix}=
V\ee{k_1}{k_2}{i_1}{i_2}V\ee{k_3}{l_3}{i_3}{j_3}V\ee{l_2}{l_1}{j_2}{j_1}$$
$$V_4\begin{pmatrix}l_1\, l_2\, l_3\, l_4\,  j_1\, j_2\, j_3\, j_4 \cr k_1 \, k_2\, k_3\, k_4\, i_1 \, i_2\, i_3\, i_4\end{pmatrix}=
V\ee{k_1}{k_2}{i_1}{i_2}V\ee{k_3}{k_4}{i_3}{i_4}
V\ee{l_4}{l_3}{j_4}{j_3}V\ee{l_2}{l_1}{j_2}{j_1}$$
$$\cdots$$

Since $\varepsilon$ and $\Delta$ are multiplicative, ($\varepsilon$) and ($\Delta$) for $V$ imply ($\varepsilon$) and ($\Delta$) for $V_2$.
\begin{eqnarray*}
\varepsilon V_2\begin{pmatrix}l_1 & l_2 & j_1 & j_2 \cr k_1 & k_2 &i_1 & i_2\end{pmatrix}
&=&\varepsilon{V}\ee{k_1}{k_2}{i_1}{i_2}\varepsilon{V}\ee{l_2}{l_1}{j_2}{j_1}\\
&=&\delta_{k_1i_1}\delta_{k_2i_2}\delta_{l_2j_2}\delta_{l_1j_1}1\\
&=&\delta_{(k_1k_2)(i_1i_2)}\delta_{(l_1l_2)(j_1j_2)}1
\end{eqnarray*}
\begin{eqnarray*}
\Delta V_2\begin{pmatrix}l_1 & l_2 & j_1 & j_2 \cr k_1 & k_2 &i_1 & i_2\end{pmatrix}
&=&\Delta V\ee{k_1}{k_2}{i_1}{i_2}\Delta V\ee{l_2}{l_1}{j_2}{j_1}\\
&=&\sum {V}\ee{k_1}{k_2}{g_1}{h_1}{V}\ee{l_2}{l_1}{g_2}{h_2}\otimes {V}\ee{g_1}{h_1}{i_1}{i_2}{V}\ee{g_2}{h_2}{j_2}{j_1}\\
&=&\sum V_2\begin{pmatrix}l_1&l_2&h_2&g_2\cr k_1&k_2&g_1&h_1\end{pmatrix}
\otimes V_2\begin{pmatrix}h_2&g_2&j_1&j_2\cr g_1&h_1&i_1&i_2\end{pmatrix}
\end{eqnarray*}

Since $S$ and $*$ are antimultiplicative, ($S$) and ($*$) for $V$ imply ($S$) and $(*)$ for $V_2$.
\begin{eqnarray*}
SV_2\begin{pmatrix}l_1 & l_2 & j_1 & j_2 \cr k_1 & k_2 &i_1 & i_2\end{pmatrix}
&=&SV\ee{l_2}{l_1}{j_2}{j_1}SV\ee{k_1}{k_2}{i_1}{i_2}\\
&=&q_{l_2}^{2}q_{j_2}^{-2}q_{j_1}^{2}q_{l_1}^{-2}
q_{k_1}^{2}q_{i_1}^{-2}q_{i_2}^{2}q_{k_2}^{-2}{V}\ee{j_1}{j_2}{l_1}{l_2}{V}\ee{i_2}{i_1}{k_2}{k_1}\\
&=&q_{(k_1k_2)}^{2}q_{(i_1i_2)}^{-2}q_{(j_1j_2)}^{2}q_{(l_1l_2)}^{-2}
V_2\begin{pmatrix}i_1&i_2&k_1&k_2\cr j_1&j_2&l_1&l_2\end{pmatrix}
\end{eqnarray*}
\begin{eqnarray*}
V_2
\begin{pmatrix}l_1 & l_2 & j_1 & j_2 \cr k_1 & k_2 &i_1 & i_2\end{pmatrix}^*
&=&V\ee{l_2}{l_1}{j_2}{j_1}^*V\ee{k_1}{k_2}{i_1}{i_2}^*\\
&=&V\ee{l_1}{l_2}{j_1}{j_2}V\ee{k_2}{k_1}{i_2}{i_1}\\
&=&V_2
\begin{pmatrix}k_1 & k_2 & i_1 & i_2 \cr l_1 & l_2 &j_1 & j_2\end{pmatrix}
\end{eqnarray*}

By using ($m^\circ$) and ($u^\circ$) for $V$ we get ($u^\circ$) for $V_2$.
\begin{eqnarray*}
\sum q_{(i_1i_2)}^2V_2
\begin{pmatrix}l_1 & l_2 & i_1 & i_2 \cr k_1 & k_2 &i_1 & i_2\end{pmatrix}
&=&\sum q_{i_1}^2q_{i_2}^{-2}V\ee{k_1}{k_2}{i_1}{i_2}V\ee{l_2}{l_1}{i_2}{i_1}\\
&=&\sum q_{i_1}^2\delta_{k_2l_2}q_{l_2}^{-2}V\ee{k_1}{l_1}{i_1}{i_1}\\
&=&\delta_{k_2l_2}q_{l_2}^{-2}\delta_{k_1l_1}q_{k_1}^{2}\\
&=&\delta_{(k_1k_2)(l_1l_2)}q_{(k_1k_2)}^2
\end{eqnarray*}

By using (${\, }^\circ m$) and (${\, }^\circ u$) for $V$ we get (${\, }^\circ m$) for $V_2$.
\begin{eqnarray*}
&&\sum q_{(s_1s_2)}^{-2}V_2\begin{pmatrix}s_1&s_2&h_1&h_2\cr k_1&k_2&g_1&g_2\end{pmatrix}
V_2\begin{pmatrix}l_1&l_2&j_1&j_2\cr s_1&s_2&i_1&i_2\end{pmatrix}\\
&=&\sum q_{s_1}^{-2}q_{s_2}^{2}
{V}\ee{k_1}{k_2}{g_1}{g_2}{V}\ee{s_2}{s_1}{h_2}{h_1}
{V}\ee{s_1}{s_2}{i_1}{i_2}{V}\ee{l_2}{l_1}{j_2}{j_1}\\
&=&\sum \delta_{h_1i_1}q_{i_1}^{-2}q_{s_2}^{2}
{V}\ee{k_1}{k_2}{g_1}{g_2}{V}\ee{s_2}{s_2}{h_2}{i_2}{V}\ee{l_2}{l_1}{j_2}{j_1}\\
&=&\delta_{h_2i_2}q_{i_2}^{2}\delta_{h_1i_1}q_{i_1}^{-2}
{V}\ee{k_1}{k_2}{g_1}{g_2}{V}\ee{l_2}{l_1}{j_2}{j_1}\\
&=&\delta_{(h_1h_2)(i_1i_2)}q_{(i_1i_2)}^{-2}V_2
\begin{pmatrix}l_1&l_2&j_1&j_2\cr k_1&k_2&g_1&g_2\end{pmatrix}
\end{eqnarray*}

Since $\varepsilon$ and $\Delta$ are multiplicative, ($\varepsilon$) and ($\Delta$) for $V$ imply ($\varepsilon$) and ($\Delta$) for $V_3$.
\begin{eqnarray*}
\varepsilon V_3\begin{pmatrix}l_1 & l_2 &l_3& j_1 & j_2&j_3 \cr k_1 & k_2&k_3 &i_1 & i_2&i_3\end{pmatrix}
&=&\varepsilon  V\ee{k_1}{k_2}{i_1}{i_2}\varepsilon V\ee{k_3}{l_3}{i_3}{j_3}\varepsilon V\ee{l_2}{l_1}{j_2}{j_1}\\
&=&\delta_{k_1i_1}\delta_{k_2i_2}\delta_{k_3i_3}\delta_{l_3j_3}\delta_{l_2j_2}\delta_{l_1j_1}1\\
&=&\delta_{(k_1k_2k_3)(i_1i_2i_3)}\delta_{(l_1l_2l_3)(j_1j_2j_3)}1
\end{eqnarray*}
\begin{eqnarray*}
&&\Delta V_3\begin{pmatrix}l_1 & l_2 &l_3& j_1 & j_2&j_3 \cr k_1 & k_2&k_3 &i_1 & i_2&i_3\end{pmatrix}\\
&=&\Delta  V\ee{k_1}{k_2}{i_1}{i_2}\Delta V\ee{k_3}{l_3}{i_3}{j_3}\Delta V\ee{l_2}{l_1}{j_2}{j_1}\\
&=&\sum {V}\ee{k_1}{k_2}{g_1}{h_1}{V}\ee{k_3}{l_3}{g_2}{h_2}V\ee{l_2}{l_1}{g_3}{h_3}\otimes {V}\ee{g_1}{h_1}{i_1}{i_2}{V}\ee{g_2}{h_2}{i_3}{j_3}V\ee{g_3}{h_3}{j_2}{j_1}\\
&=&\sum  V_3
\begin{pmatrix}l_1&l_2&l_3&h_3&g_3&h_2\cr k_1&k_2&k_3&g_1&h_1&g_2\end{pmatrix}
\otimes V_3
\begin{pmatrix}h_3&g_3&h_2&j_1&j_2&j_3\cr g_1&h_1&g_2&i_1&i_2&i_3\end{pmatrix}
\end{eqnarray*}

Since $S$ and $*$ are antimultiplicative, ($S$) and ($*$) for $V$ imply ($S$) and $(*)$ for $V_3$.
\begin{eqnarray*}
&&SV_3
\begin{pmatrix}l_1 & l_2 &l_3& j_1 & j_2&j_3 \cr k_1 & k_2&k_3 &i_1 & i_2&i_3\end{pmatrix}\\
&=&SV\ee{l_2}{l_1}{j_2}{j_1}SV\ee{k_3}{l_3}{i_3}{j_3}SV\ee{k_1}{k_2}{i_1}{i_2}\\
&=&q_{l_2}^{2}q_{j_2}^{-2}q_{j_1}^{2}q_{l_1}^{-2}
q_{k_3}^{2}q_{i_3}^{-2}q_{j_3}^{2}q_{l_3}^{-2}
q_{k_1}^{2}q_{i_1}^{-2}q_{i_2}^{2}q_{k_2}^{-2}
{V}\ee{j_1}{j_2}{l_1}{l_2}V\ee{j_3}{i_3}{l_3}{k_3}{V}\ee{i_2}{i_1}{k_2}{k_1}\\
&=&q_{(k_1k_2k_3)}^{2}q_{(i_1i_2i_3)}^{-2}q_{(j_1j_2j_3)}^{2}q_{(l_1l_2l_3)}^{-2}
V_3\begin{pmatrix}i_1&i_2&i_3&k_1&k_2&k_3\cr j_1&j_2&j_3&l_1&l_2&l_3\end{pmatrix}
\end{eqnarray*}
\begin{eqnarray*}
V_3\begin{pmatrix}l_1 & l_2 &l_3& j_1 & j_2&j_3 \cr k_1 & k_2&k_3 &i_1 & i_2&i_3\end{pmatrix}^*
&=&V\ee{l_2}{l_1}{j_2}{j_1}^*V\ee{k_3}{l_3}{i_3}{j_3}^*V\ee{k_1}{k_2}{i_1}{i_2}^*\\
&=&V\ee{l_1}{l_2}{j_1}{j_2}V\ee{l_3}{k_3}{j_3}{i_3}V\ee{k_2}{k_1}{i_2}{i_1}\\
&=&V_3\begin{pmatrix}k_1 & k_2&k_3& i_1 & i_2&i_3\cr l_1 & l_2 &l_3 &j_1 & j_2&j_3 \end{pmatrix}
\end{eqnarray*}

By using ($u^\circ$), ($m^\circ$) and ($u^\circ$) again for $V$ we get ($u^\circ$) for $V_3$.
\begin{eqnarray*}
&&\sum q_{(i_1i_2i_3)}^2V_3
\begin{pmatrix}l_1 & l_2&l_3 & i_1 & i_2&i_3 \cr k_1 & k_2&k_3 &i_1 & i_2&i_3\end{pmatrix}\\
&=&\sum q_{i_1}^2q_{i_2}^{-2}q_{i_3}^2V\ee{k_1}{k_2}{i_1}{i_2}V\ee{k_3}{l_3}{i_3}{i_3}V\ee{l_2}{l_1}{i_2}{i_1}\\
&=&\sum q_{i_1}^2q_{i_2}^{-2}\delta_{k_3l_3}q_{k_3}^2V\ee{k_1}{k_2}{i_1}{i_2}V\ee{l_2}{l_1}{i_2}{i_1}\\
&=&\sum q_{i_1}^2\delta_{k_3l_3}q_{k_3}^2\delta_{k_2l_2}q_{l_2}^{-2}V\ee{k_1}{l_1}{i_1}{i_1}\\
&=&\delta_{k_3l_3}q_{k_3}^2\delta_{k_2l_2}q_{l_2}^{-2}\delta_{k_1l_1}q_{k_1}^{2}\\
&=&\delta_{(k_1k_2k_3)(l_1l_2l_3)}q_{(k_1k_2k_3)}^2
\end{eqnarray*}

By using (${\, }^\circ m$), (${\, }^\circ u$) and (${\, }^\circ m$) again for $V$ we get (${\, }^\circ m$) for $V_3$.
\begin{eqnarray*}
&&\sum q_{(s_1s_2s_3)}^{-2}V_3\begin{pmatrix}s_1&s_2&s_3&h_1&h_2&h_3\cr k_1&k_2&k_3&g_1&g_2&g_3\end{pmatrix}
V_3\begin{pmatrix}l_1&l_2&l_3&j_1&j_2&j_3\cr s_1&s_2&s_3&i_1&i_2&i_3\end{pmatrix}\\
&=&\sum q_{s_1}^{-2}q_{s_2}^{2}q_{s_3}^{-2}\\
&&\cdot\,\, {V}\ee{k_1}{k_2}{g_1}{g_2}
{V}\ee{k_3}{s_3}{g_3}{h_3}
{V}\ee{s_2}{s_1}{h_2}{h_1}
{V}\ee{s_1}{s_2}{i_1}{i_2}
{V}\ee{s_3}{l_3}{i_3}{j_3}
{V}\ee{l_2}{l_1}{j_2}{j_1}\\
&=&\sum q_{s_2}^{2}q_{s_3}^{-2}\delta_{h_1i_1}q_{i_1}^{-2}
{V}\ee{k_1}{k_2}{g_1}{g_2}
{V}\ee{k_3}{s_3}{g_3}{h_3}
{V}\ee{s_2}{s_2}{h_2}{i_2}
{V}\ee{s_3}{l_3}{i_3}{j_3}
{V}\ee{l_2}{l_1}{j_2}{j_1}\\
&=&\sum q_{s_3}^{-2}\delta_{h_1i_1}q_{i_1}^{-2}\delta_{h_2i_2}q_{i_2}^{2}
{V}\ee{k_1}{k_2}{g_1}{g_2}
{V}\ee{k_3}{s_3}{g_3}{h_3}
{V}\ee{s_3}{l_3}{i_3}{j_3}
{V}\ee{l_2}{l_1}{j_2}{j_1}\\
&=&\delta_{h_1i_1}q_{i_1}^{-2}\delta_{h_2i_2}q_{i_2}^{2}\delta_{h_3i_3}q_{i_3}^{-2}
{V}\ee{k_1}{k_2}{g_1}{g_2}
{V}\ee{k_3}{l_3}{g_3}{j_3}
{V}\ee{l_2}{l_1}{j_2}{j_1}\\
&=&\delta_{(h_1h_2h_3)(i_1i_2i_3)}q_{(i_1i_2i_3)}^{-2}V_3
\begin{pmatrix}l_1&l_2&l_3&j_1&j_2&j_3\cr k_1&k_2&k_3&g_1&g_2&g_3\end{pmatrix}
\end{eqnarray*}

The proof for arbitrary $n$ even is similar to the proof for $n=2$ and for arbitrary $n$ odd, to the proof for $n=3$.
\end{proof}

A linear map $T:A^{\otimes n}\to A^{\otimes m}$ is $v_\infty$-equivariant if the following diagram commutes.
$$\begin{matrix}A^{\otimes n}& \displaystyle{\mathop{\longrightarrow}^{T}}&
  A^{\otimes m}\cr &\ \cr v_n\downarrow & \ & \downarrow v_m\cr&\ \cr 
A^{\otimes n}\otimes H& \displaystyle{\mathop{\longrightarrow}^{T\otimes id}}&
 A^{\otimes m}\otimes H \end{matrix}$$

For $n=0$ a map $T$ is $v_\infty$-equivariant if and only if $T(1)$ is fixed by $v_m$.

\begin{lemm}
The following linear maps
$${I_n}\begin{pmatrix}j_1 & \ldots & j_{n-1}\cr i_1 & \ldots & i_{n-1}\end{pmatrix}=
\sum\begin{pmatrix}j_1 & \ldots & j_{n-1} & l\cr i_1 & \ldots & i_{n-1}& l\end{pmatrix}$$
$$\wt{e}_n=\sum q_i^{\pm 2}q_j^{\pm 2}\begin{pmatrix}g_1 & \ldots & g_{n-2} & j & j\cr g_1 & \ldots &
  g_{n-2} & i & i\end{pmatrix}$$
$$\wt{E}_n\begin{pmatrix}j_1 & \ldots & j_n\cr i_1 & \ldots & i_n\end{pmatrix}=
\delta_{i_nj_n}q_{i_n}^{\mp 4}\begin{pmatrix}j_1 & \ldots & j_{n-1}\cr i_1 & \ldots & i_{n-1}\end{pmatrix}$$
where $\pm 1=(-1)^n$, are $v_\infty$-equivariant.
\end{lemm}

\begin{proof}
The coactions $v_2$ and $v_3$ are given by the following formulae.
\begin{eqnarray*}
v_2\begin{pmatrix}j_1&j_2\cr i_1&i_2\end{pmatrix}
&=&\sum \begin{pmatrix}l_1&l_2\cr k_1&k_2\end{pmatrix}
\otimes q_{k_1}^{-1}q_{k_2}q_{i_1}q_{i_2}^{-1}q_{j_1}q_{j_2}^{-1}
q_{l_1}^{-1}q_{l_2}V\ee{k_1}{k_2}{i_1}{i_2}V\ee{l_2}{l_1}{j_2}{j_1}\\
v_3\begin{pmatrix}j_1&j_2&j_3\cr i_1&i_2&i_3\end{pmatrix}
&=&\sum \begin{pmatrix}l_1&l_2&l_3\cr k_1&k_2&k_3\end{pmatrix}
\otimes V\ee{k_1}{k_2}{i_1}{i_2}V\ee{k_3}{l_3}{i_3}{j_3}V\ee{l_2}{l_1}{j_2}{j_1}\\
&&\cdot\,\, 
q_{k_1}^{-1}q_{k_2}q_{k_3}^{-1}q_{i_1}q_{i_2}^{-1}q_{i_3}q_{j_1}q_{j_2}^{-1}q_{j_3}
q_{l_1}^{-1}q_{l_2}q_{l_3}^{-1}
\end{eqnarray*}

By using ($m^\circ$) we get that ${I_2}$ is $v_\infty$-equivariant.
\begin{eqnarray*}
v_2{I_2}\e{i_1}{j_1}
&=&\sum\ee{k_1}{l_1}{k_2}{l_2}\otimes
q_{k_1}^{-1}q_{k_2}q_{i_1}q_{j_1}
q_{l_1}^{-1}q_{l_2}q_{l}^{-2}{V}\ee{k_1}{k_2}{i_1}{l}{V}\ee{l_2}{l_1}{l}{j_1}\\
&=&\sum\ee{k_1}{l_1}{k_2}{l_2}\otimes
\delta_{k_2l_2}q_{k_2}^{-2}q_{k_1}^{-1}q_{k_2}q_{i_1}q_{j_1}
q_{l_1}^{-1}q_{l_2}{V}\ee{k_1}{l_1}{i_1}{j_1}\\
&=&\sum \ee{k_1}{l_1}{l}{l}\otimes 
q_{k_1}^{-1}q_{i_1}q_{j_1}q_{l_1}^{-1}
{V}\ee{k_1}{l_1}{i_1}{j_1}\\
&=&({I_2}\otimes id)v\e{i_1}{j_1}
\end{eqnarray*}

By using ($u^\circ$) we get that ${I_3}$ is $v_\infty$-equivariant.
\begin{eqnarray*}
v_3{I_3}\ee{i_1}{j_1}{i_2}{j_2}
&=&\sum\eee{k_1}{l_1}{k_2}{l_2}{k_3}{l_3}\otimes
{V}\ee{k_1}{k_2}{i_1}{i_2}
{V}\ee{k_3}{l_3}{l}{l}
V\ee{l_2}{l_1}{j_2}{j_1}\\
&&\cdot\,\, 
q_{k_1}^{-1}q_{k_2}q_{k_3}^{-1}q_{i_1}q_{i_2}^{-1}q_{j_1}q_{j_2}^{-1}
q_{l_1}^{-1}q_{l_2}q_{l_3}^{-1}q_l^2\\
&=&\sum\eee{k_1}{l_1}{k_2}{l_2}{k_3}{l_3}\otimes{V}\ee{k_1}{k_2}{i_1}{i_2}
V\ee{l_2}{l_1}{j_2}{j_1}\\
&&\cdot\,\, \delta_{l_3k_3}
q_{k_3}^2q_{k_1}^{-1}q_{k_2}q_{k_3}^{-1}q_{i_1}q_{i_2}^{-1}q_{j_1}q_{j_2}^{-1}
q_{l_1}^{-1}q_{l_2}q_{l_3}^{-1}q_l^2\\
&=&\sum\eee{k_1}{l_1}{k_2}{l_2}{l}{l}\otimes{V}\ee{k_1}{k_2}{i_1}{i_2}
V\ee{l_2}{l_1}{j_2}{j_1}\\
&&\cdot\,\, q_{k_1}^{-1}q_{k_2}q_{i_1}q_{i_2}^{-1}q_{j_1}q_{j_2}^{-1}
q_{l_1}^{-1}q_{l_2}q_l^2\\
&=&({I_3}\otimes id)v_2\ee{i_1}{j_1}{i_2}{j_2}
\end{eqnarray*}

By using ($u^\circ$) twice we get that $\wt{e}_2$ is $v_\infty$-equivariant.
\begin{eqnarray*}
v_2(\wt{e}_2)
&=&\sum\ee{k_1}{l_1}{k_2}{l_2}\otimes
q_i^2q_j^2q_{k_1}^{-1}q_{k_2}q_{l_1}^{-1}q_{l_2}
{V}\ee{k_1}{k_2}{i}{i}{V}\ee{l_2}{l_1}{j}{j}\\
&=&\sum \ee{k_1}{l_1}{k_2}{l_2}\otimes \delta_{k_1k_2}\delta_{l_1l_2}q_{k_1}^2q_{l_1}^2\\
&=&\wt{e}_2\otimes 1
\end{eqnarray*}

By using ($m^\circ$) twice and ($u^\circ$) we get that $\wt{e}_3$ is $v_\infty$-equivariant.
\begin{eqnarray*}
v_3(\wt{e}_3)
&=&\sum\eee{k_1}{l_1}{k_2}{l_2}{k_3}{l_3}\otimes
{V}\ee{k_1}{k_2}{g_1}{i}
{V}\ee{k_3}{l_3}{i}{j}
V\ee{l_2}{l_1}{j}{g_1}\\
&&\cdot\,\, q_i^{-2}q_j^{-2}q_{k_1}^{-1}q_{k_2}q_{k_3}^{-1}
q_{g_1}^2q_{l_1}^{-1}q_{l_2}q_{l_3}^{-1}\\
&=&\sum\eee{k_1}{l_1}{k_2}{l_2}{k_3}{l_3}\otimes \delta_{l_2l_3}
q_i^{-2}q_{l_2}^{-2}q_{k_1}^{-1}q_{k_2}q_{k_3}^{-1}
q_{g_1}^2q_{l_1}^{-1}
{V}\ee{k_1}{k_2}{g_1}{i}
{V}\ee{k_3}{l_1}{i}{g_1}\\
&=&\sum\eee{k_1}{l_1}{k_2}{l_2}{k_3}{l_3}\otimes\delta_{l_2l_3}\delta_{k_2k_3}
q_{k_2}^{-2}q_{l_2}^{-2}q_{k_1}^{-1}q_{g_1}^2q_{l_1}^{-1}
{V}\ee{k_1}{l_1}{g_1}{g_1}\\
&=&\sum\eee{k_1}{l_1}{k_2}{l_2}{k_3}{l_3}\otimes
\delta_{k_1l_1}\delta_{l_2l_3}\delta_{k_2k_3}q_{k_2}^{-2}q_{l_2}^{-2}\\
&=&\wt{e}_3\otimes 1
\end{eqnarray*}

By using (${\, }^\circ m$) we get that $\wt{E}_2$ is $v_\infty$-equivariant.
\begin{eqnarray*}
(\wt{E}_2\otimes id)v_2\ee{i_1}{j_1}{i_2}{j_2}
&=&\sum\e{k_1}{l_1}\otimes q_g^{-2}q_{k_1}^{-1}q_{i_1}q_{i_2}^{-1}q_{j_1}q_{j_2}^{-1}
q_{l_1}^{-1}
{V}\ee{k_1}{g}{i_1}{i_2}{V}\ee{g}{l_1}{j_2}{j_1}\\
&=&\sum\e{k_1}{l_1}\otimes
\delta_{i_2j_2} q_{i_2}^{-4}q_{k_1}^{-1}q_{i_1}q_{j_1}q_{l_1}^{-1}
{V}\ee{k_1}{l_1}{i_1}{j_1}\\
&=&v \wt{E}_2\ee{i_1}{j_1}{i_2}{j_2}
\end{eqnarray*}

By using (${\, }^\circ u$) we get that $\wt{E}_3$ is $v_\infty$-equivariant.
\begin{eqnarray*}
(\wt{E}_3\otimes id)v_3\eee{i_1}{j_1}{i_2}{j_2}{i_3}{j_3}
&=&\sum\ee{k_1}{l_1}{k_2}{l_2}\otimes
{V}\ee{k_1}{k_2}{i_1}{i_2}
{V}\ee{g}{g}{i_3}{j_3}
V\ee{l_2}{l_1}{j_2}{j_1}\\
&&\cdot\,\, 
q_g^{2}q_{k_1}^{-1}q_{k_2}q_{i_1}q_{i_2}^{-1}q_{i_3}q_{j_1}q_{j_2}^{-1}q_{j_3}
q_{l_1}^{-1}q_{l_2}\\
&=&\sum\ee{k_1}{l_1}{k_2}{l_2}\otimes {V}\ee{k_1}{k_2}{i_1}{i_2}
V\ee{l_2}{l_1}{j_2}{j_1}\\
&&\cdot\,\, \delta_{i_3j_3}
q_{i_3}^4q_{k_1}^{-1}q_{k_2}q_{i_1}q_{i_2}^{-1}q_{j_1}q_{j_2}^{-1}
q_{l_1}^{-1}q_{l_2}\\
&=&v_2 \wt{E}_3\eee{i_1}{j_1}{i_2}{j_2}{i_3}{j_3}
\end{eqnarray*}

The proof for arbitrary $n$ even is similar to the proof for $n=2$ and for arbitrary $n$ odd, to the proof for $n=3$.
\end{proof}

Let $h:H\to\c$ be the Haar integral constructed by Woronowicz in \cite{w1}. This is a unital linear form having the following ``bi-invariance'' property.
$$({h}\otimes id)\Delta =(id\otimes{h} )\Delta ={h} (.)1$$

If $v :A\to A\otimes H$ is a coaction then $\Gamma_n =(id\otimes h)v_n$ is an idempotent of $\l (A^{\otimes n})$ and its image are the fixed points of $v_n$. This follows from the computation
$$v_n \Gamma_n(x)=v_n (id\otimes{h} )v_n (x)
=(id\otimes id\otimes{h} )(id\otimes\Delta )v_n (x)=\Gamma_n(x)\otimes 1$$

A pair of linear maps $(T,T^q): A^{\otimes n}\to A^{\otimes m}$ is called ``weakly $v_\infty$-equivariant'' if the following diagram commutes.
$$\begin{matrix}A^{\otimes n}& \displaystyle{\mathop{\longrightarrow}^{T^q}}&
  A^{\otimes m}\cr &\ \cr \Gamma_{n}\downarrow & \ & \downarrow \Gamma_{m}\cr &\ \cr
A^{\otimes n}& \displaystyle{\mathop{\longrightarrow}^{T}}&
 A^{\otimes m}\end{matrix}$$

The interest in this notion is that it makes the following diagram factorise.
$$\begin{matrix}A^{\otimes n}& \displaystyle{\mathop{\longrightarrow}^{T}} &
 A^{\otimes m}\cr &\ \cr
\cup & \ & \cup\cr &\ \cr
Im(\Gamma_{n}) & \displaystyle{\mathop{\longrightarrow}}&
Im(\Gamma_{m})\end{matrix}$$

We say that an operator $T$ is ``weakly $v_\infty$-equivariant'' if the pair $(T,T)$ is weakly $v_\infty$-equivariant. This happens for instance if $T$ is $v_\infty$-equivariant, because we can glue the $v_\infty$-equivariance diagram of $T$ to the following trivial diagram.
$$\begin{matrix}A^{\otimes n}\otimes H& \displaystyle{\mathop{\longrightarrow}^{T\otimes id}} &
 A^{\otimes m}\otimes H\cr &\ \cr
id\otimes h\downarrow & \ & \downarrow id\otimes h\cr &\ \cr
A^{\otimes n}& \displaystyle{\mathop{\longrightarrow}^{T}}&
 A^{\otimes m}\end{matrix}$$

The following linear map, called modular map of $\varphi$
$$\theta\e{i}{j}=q_i^4q_j^{-4}\e{i}{j}$$
is the unique linear map $\theta :A\to A$ such that $\varphi (ab)=\varphi (b\theta (a))$ for any $a,b\in A$.

Consider the automorphism $\sigma :H\to H$ constructed by Woronowicz in \cite{w1}, which satisfies $h(ab)=h(b\sigma (a))$ for any $a,b$.

\begin{lemm}
Assume that the following ``modularity'' condition is satisfied.
$$(\theta\otimes id)v\theta =(id\otimes\sigma )v$$

(i) The following pair of maps is weakly $v_\infty$-equivariant.
$$J_n\begin{pmatrix}j_3 & \ldots & j_{n+1}\cr i_3 & \ldots & i_{n+1}\end{pmatrix}=
\sum\begin{pmatrix}l&k&j_3&\ldots&j_{n+1}\cr l&k&i_3 & \ldots & i_{n+1}\end{pmatrix}$$
$$J_n^q\begin{pmatrix}j_3 & \ldots & j_{n+1}\cr i_3 & \ldots & i_{n+1}\end{pmatrix}=
\sum q_l^{-8}q_k^8
\begin{pmatrix}l&k&j_3&\ldots&j_{n+1}\cr l&k&i_3 & \ldots & i_{n+1}\end{pmatrix}$$

(ii) If $\varphi$ has the trace property $\varphi (ab)=\varphi (ba)$ then $J_n$ is weakly $v_\infty$-equivariant.

(iii) If $H$ is commutative then $J_n$ is $v_\infty$-equivariant.
\end{lemm}

\begin{proof}
(i) By using the formulae of $v$ and $\theta$ we get
$$(\theta\otimes id)v\theta \e{i}{j}=\sum \e{k}{l}\otimes 
q_k^4q_l^{-4}q_i^4q_j^{-4}\cdot q_k^{-1}q_iq_jq_l^{-1}V\ee{k}{l}{i}{j}$$
so the modularity condition is equivalent to the following condition ($\sigma$).
$$\sigma V\ee{k}{l}{i}{j}=q_k^4q_l^{-4}q_i^4q_j^{-4} V\ee{k}{l}{i}{j}$$

By using the formula of $J_n^q$ we get
$$v_{n+1}J_n^q
\begin{pmatrix}j_3& \ldots &j_{n+1}\cr i_3&\ldots& i_{n+1}\end{pmatrix}
=v_{n+1}\sum q_l^{-8}q_k^8\begin{pmatrix}l&k&j_3&\ldots&j_{n+1}\cr l&k&i_3 & \ldots & i_{n+1}\end{pmatrix}$$
$$=\sum \begin{pmatrix}l_1&l_2&l_3&\ldots&l_{n+1}\cr k_1&k_2&k_3& \ldots & k_{n+1}\end{pmatrix}\otimes q_{(k_3\ldots k_{n+1})}^{-1}q_{(i_3\ldots i_{n+1})}
q_{(j_3\ldots j_{n+1})}q_{(l_3\ldots l_{n+1})}^{-1}Z$$
with $Z$ given by the following formula.
$$Z=\sum q_l^{-6}q_k^6q_{k_1}^{-1}q_{k_2}q_{l_1}^{-1}q_{l_2}V\ee{k_1}{k_2}{l}{k}
V\ee{k_3}{k_4}{i_3}{i_4}\ldots  V\ee{l_4}{l_3}{j_4}{j_3}
V\ee{l_2}{l_1}{k}{l}$$

By applying the Haar integral to $Z$ we get
$$h(Z)=h\left( V\ee{k_3}{k_4}{i_3}{i_4}\ldots  V\ee{l_4}{l_3}{j_4}{j_3}T\right)$$
with $T$ given by the following formula.
$$T=\sum q_l^{-6}q_k^6q_{k_1}^{-1}q_{k_2}q_{l_1}^{-1}q_{l_2}V\ee{l_2}{l_1}{k}{l}\sigma V\ee{k_1}{k_2}{l}{k}$$

By using ($\sigma$), ($m^\circ$) and ($u^\circ$) we can compute $T$.
\begin{eqnarray*}
T&=&\sum q_l^{-6}q_k^6q_{k_1}^{-1}q_{k_2}q_{l_1}^{-1}q_{l_2}\cdot q_{k_1}^4q_l^{4}q_k^{-4}q_{k_2}^{-4}V\ee{l_2}{l_1}{k}{l}V\ee{k_1}{k_2}{l}{k}\\
&=&\sum q_l^{-2}q_k^2q_{k_1}^{3}q_{k_2}^{-3}q_{l_1}^{-1}q_{l_2}V\ee{l_2}{l_1}{k}{l}V\ee{k_1}{k_2}{l}{k}\\
&=&\sum \delta_{k_1l_1}q_{l_1}^{-2}q_k^2q_{k_1}^{3}q_{k_2}^{-3}q_{l_1}^{-1}q_{l_2}
V\ee{k_2}{l_2}{k}{k}\\
&=&\delta_{k_1l_1}q_{l_1}^{-2}\delta_{k_2l_2}q_{k_2}^{2}q_{k_1}^{3}q_{k_2}^{-3}q_{l_1}^{-1}q_{l_2}\\
&=&\delta_{k_1l_1}\delta_{k_2l_2}
\end{eqnarray*}

Thus by applying $id\otimes h$ to the formula of $v_{n+1}J_n^q$ we get
\begin{eqnarray*}
&&(id\otimes h)v_{n+1}J_n^q
\begin{pmatrix}j_3& \ldots &j_{n+1}\cr i_3&\ldots& i_{n+1}\end{pmatrix}\\
&=&\sum \begin{pmatrix}k_1&l_2&l_3&\ldots&l_{n+1}\cr k_1&l_2&k_3& \ldots & k_{n+1}\end{pmatrix}
h\left( V\ee{k_3}{k_4}{i_3}{i_4}\ldots  V\ee{l_4}{l_3}{j_4}{j_3}\right)\\
&&\cdot\,\, q_{(k_3\ldots k_{n+1})}^{-1}q_{(i_3\ldots i_{n+1})}
q_{(j_3\ldots j_{n+1})}q_{(l_3\ldots l_{n+1})}^{-1}\\
&=&(id\otimes h)({J_n}\otimes id)v_{n-1}
\begin{pmatrix}j_3& \ldots &j_{n+1}\cr i_3&\ldots& i_{n+1}\end{pmatrix}
\end{eqnarray*}
so the weak $v_\infty$-equivariance diagram commutes.

(ii) In terms of weights, the fact that $\varphi$ is a trace means that $q_i^4$ depends only on the matrix block containing $(^i_i)$. Thus the spin factor $q_l^{-8}q_k^8$ in the formula of $J_n^q$ cancels, because $(^l_l)$ and $(^k_k)$ are in the same matrix block of $A$ (cf. loop notation).

(iii) The Haar integral and its modular map are used in proof of (i) for rotating a product of elements of $H$. If $H$ is commutative its product does the same job.
\end{proof}

We have all ingredients needed for the case of trace-preserving coactions. Let $\p$ be the colored planar operad constructed by Jones in \cite{j1}. A planar algebra is a sequence of vector spaces $P=P_0^\pm ,\, P_1,\, P_2,\,\, P_3,\ldots$ with a colored operad morphism $\pi : \p\to Hom(P)$, where $Hom(P)$ is the colored operad of multilinear maps between $P_n$'s.

If $Q_n\subset P_n$ is a sequence of subspaces, the restriction of multiplinear maps between $P_n$'s to multilinear maps between $Q_n$'s is a partially defined colored operad morphism $Res: Hom(P)\to Hom(Q)$. If the domain of $Res$ contains the image of $\pi$ the composition of $Res$ and $\pi$ makes $Q$ a planar algebra, called subalgebra of $P$.

The annular category $A$ is defined as follows. The objects are the positive integers and the space $A(i,j)$ of arrows from $i$ to $j$ is formed by tangles in $\p$ with ``output'' disc having $2j$ marked points and one ``input'' disc, having $2i$ marked points. Composition of arrows is given by gluing of annuli. The restriction of $\pi$ to $A$ is a morphism from $A$ to the category $\l (P)$ having as arrows from $i$ to $j$ the linear maps from $P_i$ to $P_j$.

Let $(A,\varphi )$ be as in \S 1 and assume that $\varphi$ has the trace property $\varphi (ab)=\varphi (ba)$. Let $\pi :\p\to Hom(P(A,\varphi ))$ be the planar algebra associated to the bipartite graph of $A$, with spin vector $a\mapsto\varphi (1_a)$. The sequence of vector spaces of $P(A,\varphi )$ will be canonically identified with the sequence of tensor powers of $A$. See Jones \cite{j2}.

\begin{theo}
(i) The linear maps in the image by $\pi$ of the annular category are weakly $v_\infty$-equivariant. If $H$ is commutative, they are $v_\infty$-equivariant.

(ii) The spaces of fixed points of the coactions $v_n$ form a subalgebra of $P(A,\varphi )$.
\end{theo}

\begin{proof}
(i) Gluing of commutative diagrams shows that weak $v_\infty$-equivariance is stable by composition, so the annular tangles whose image by $\pi$ are weakly $v_\infty$-equivariant from a subcategory $B\subset A$. We want to prove that $B=A$.

Consider the inclusion tangle in $A(n-1,n)$, expectation tangle in $A(n,n-1)$, Jones projection tangle in $A(0,n)$ and shift tangle in $A(n-1,n+1)$ (see \cite{j1} for pictures). By \cite{j2} their images by $\pi$ are given by the formulae in lemma 2.1 and lemma 2.2, suitably rescaled. Lemma 2.1 shows that the first three tangles are in $B$. In terms of weights, the fact that $\varphi$ is a trace means that $q_i$ depends only on the matrix block containing $(^i_i)$, so the spin factor in the formula ($S$) of the antipode cancels. In particular the square of the antipode is the identity on the coefficients of $v$, and by replacing $H$ with its $*$-subalgebra generated by these coefficients we may assume that $S^2=id$. By \cite{w1} the Haar integral is a trace, so the modularity condition in lemma 2.2 is satisfied. Thus the shift tangle is in $B$.

The sets $A(0,n)$ of Temperley-Lieb tangles being generated by inclusions and Jones projections, they are in $B$. For $x,y\in A(0,n)$ let $M(x,y)\in A(n,n)$ be the 3-multiplication $n$-tangle of $\p$ with the upper circle filled with $y$ and the lower circle filled with $x$. The corresponding linear map is $\wt{M}(x,y):p\mapsto \wt{x}p\wt{y}$ and since fixed points of $v_n$ are stable under multiplication, $M(x,y)$ is in $B$.

Let $T\in A(i,j)$. By using boxes instead of discs, as in \cite{j1}, isotope $T$, then cut it horizontally in three parts such that the middle part contains the inner box plus vertical strings only. By adding contractible circles at right, we can arrange such that the number of points on the middle cuts is greater than $j$. By adding more contractible circles at right, each of them consisting of ``up'' and ``down'' semicircles plus two outside ``expectation'' strings connecting them, we get an equality of the form $T^{\circ\circ\ldots\circ} =EM(x,y)IJ$, where $I$ is a composition of inclusion tangles, $J$ is a composition of shift tangles, $E$ is a composition of expectation tangles, $x$ and $y$ are in $A(k,k)$ for some big $k$ and $T^{\circ\circ\ldots\circ}$ is obtained from $T$ by adding contractible circles. Thus $T^{\circ\circ\ldots\circ}$ is in $B$, so $T$ is in $B$. Same proof works for the second part, with ``weak $v_\infty$-equivariant'' replaced by ``$v_\infty$-equivariant''.

(ii) Since weakly $v_\infty$-equivariant maps send fixed points to fixed points, part (i) shows that $\pi (A)$ is in the domain of $Res$. By proposition 1.18 in \cite{j1} this implies that $\pi (\p)$ is in the domain of $Res$ and this proves the first assertion.
\end{proof}

\section{Twisted structure}

We assume that $\varphi$ is a $\delta$-form, in the sense that $\varphi (1)=1$ and $Tr(B^{-4})=\delta^2$ for any matrix block $B$ of the unique $Q$ such that $\varphi =Tr(Q^4.)$. See \cite{fc} for examples and comments. In terms of the basis, the unitality of $\varphi$ translates into the following formula $(\ddag)$.
$$\sum q_j^{4}=1$$

We use the equivalence relation $i\sim j$ if $(^i_i)$ and $(^j_j)$ are in the same matrix block of $A$. For any $i$ in the set of indices we have the following formula $(\dag)$.
$$\sum_{j\sim i} q_j^{-4}=\delta^2$$

By using $\delta$ we define normalised forms, expectations and Jones projections by
$$\varphi_n\begin{pmatrix}j_1 & \ldots & j_n\cr i_1 & \ldots & i_n\end{pmatrix}
=\delta^{\frac{1}{2}\mp\frac{1}{2}-n}\delta_{(i_1\ldots i_n)(j_1\ldots j_n)}q_{(i_1\ldots i_n)}^4$$ 
$$E_n\begin{pmatrix}j_1 & \ldots & j_n\cr i_1 & \ldots & i_n\end{pmatrix}=
\delta_{i_nj_n}\delta^{-1\mp 1}q_{i_n}^{\mp 4}\begin{pmatrix}j_1 & \ldots & j_{n-1}\cr i_1 & \ldots & i_{n-1}\end{pmatrix}$$
$$e_n=\sum \delta^{-1\pm 1}q_i^{\pm 2}q_j^{\pm 2}\begin{pmatrix}g_1 & \ldots & g_{n-2} & j & j\cr g_1 & \ldots & g_{n-2} & i & i\end{pmatrix}$$
where $\pm 1=(-1)^n$. Define also a linear form $\psi_2$ by
$$\psi_2\begin{pmatrix}j_1&j_2\cr i_1&i_2\end{pmatrix}=\delta_{(i_1i_2)(j_1j_2)}
\delta^{-2}q_{i_1}^{-4}q_{i_2}^{4}$$

The modular map $\theta_n$ of $\wt{\varphi}_n$ is given by the following formula (see \S 2).
$$\theta_n\begin{pmatrix}j_1 & \ldots & j_n\cr i_1 & \ldots & i_n\end{pmatrix}
=q_{(i_1\ldots i_n)}^4q_{(j_1\ldots j_n)}^{-4}
\begin{pmatrix}j_1 & \ldots & j_n\cr i_1 & \ldots & i_n\end{pmatrix}$$

In this section we prove the following technical result.

\begin{prop}
If $Q_n\subset A^{\otimes n}$ is a sequence of $\c^*$-algebras satisfying

(1) $I_n(Q_{n-1})\subset Q_n$, $E_n(Q_{n})\subset Q_{n-1}$, $J_n(Q_{n-1})\subset Q_{n+1}$ and $e_n\in Q_n$ for any $n$.

(2) $\theta_n (x)=x$ for any $x\in Q_n$ and any $n$.

(3) $\varphi_2 (x)=\psi_2(x)$ for any $x\in Q_2$.

\noindent then there exists a unique $\c^*$-planar algebra structure on the sequence $Q_n$ such that the inclusions, shifts, traces, expectations and Jones projections are the restrictions of $I_n$, $J_n$, $\varphi_n$, $E_n$ and $e_n$. This $\c^*$-planar algebra is spherical and of modulus $\delta$.
\end{prop}

This will be proved by using the ``bubbling'' result of Jones in \cite{j1} applied to a certain lattice of $\c^*$-algebras satisfying the axioms of Popa in \cite{po}. For checking the axioms we have to verify all relevant formulae satisfied by $I_n$, $J_n$, $\varphi_n$, $E_n$ and $e_n$. This kind of computation appears in many places in the subfactor literature, see e.g. the books \cite{ghj} or \cite{js}.

In this paper the set of parameters is somehow maximal, so we will give self-contained complete proofs for everything. It is possible to use some Hopf algebra dualities in order to cut from computations, but this rather complicates things and we prefer to use the obvious symmetries only. In fact, the interesting thing would be to have an explicit construction of the partition function, as in \cite{j2} and we don't know if this is possible.

Note also that all formulae to be verified are irrelevant once the result is proved, because they can be easily verified on pictures.

We first associate to $A$ and to the numbers $q_i$ a system of $\c^*$-algebras satisfying some of Popa's axioms. Define linear maps
$$J_n^-:A^{\otimes n-1}\to A\otimes A^{\otimes n-1}\hskip 20mm 
J_n^+:A\otimes A^{\otimes n-1}\to A^{\otimes n+1}$$
by the following formulae in terms of the basis.
$$J_n^-\begin{pmatrix}j_3 & \ldots & j_{n+1}\cr i_3 & \ldots & i_{n+1}\end{pmatrix}= 
\sum\e{g}{g}\otimes \begin{pmatrix}j_3 & \ldots & j_{n+1}\cr i_3 & \ldots & i_{n+1}\end{pmatrix}$$ 
$$J_n^+\left( 
\e{i_2}{j_2}\otimes \begin{pmatrix}j_3 & \ldots & j_{n+1}\cr i_3 & \ldots & i_{n+1}\end{pmatrix}\right) =
\sum\begin{pmatrix}h&j_2& \ldots & j_{n+1}\cr h&i_2& \ldots & i_{n+1}\end{pmatrix}$$

These are inclusions of $\c^*$-algebras.

\begin{lemm}
We have $J_n^+J_n^-=J_n$ and the following diagram (I) commutes
$$\begin{matrix}
\c &\displaystyle{\mathop{\longrightarrow}^{I_1}} & A & 
\displaystyle{\mathop{\longrightarrow}^{I_2}} & A^{\otimes 2}&
\displaystyle{\mathop{\longrightarrow}^{I_3}} & A^{\otimes 3}& 
\displaystyle{\mathop{\longrightarrow}^{I_4}} & A^{\otimes 4}&\cdots\cr\ \cr
& & \uparrow J_0^+& & \uparrow J_1^+ & & \uparrow J_2^+& & \uparrow J_3^+\cr\ \cr
\ &\ & \c & 
\displaystyle{\mathop{\longrightarrow}^{id\otimes I_0}} & A\otimes\c&
\displaystyle{\mathop{\longrightarrow}^{id\otimes I_1}} & A\otimes A& 
\displaystyle{\mathop{\longrightarrow}^{id\otimes I_2}} & A\otimes A^{\otimes 2}&\cdots\cr\ \cr
& & \ & & \uparrow J_1^- & & \uparrow J_2^-& & \uparrow J_3^-\cr\ \cr
\ &\ & \ & \ & \c&
\displaystyle{\mathop{\longrightarrow}^{I_1}} & A& 
\displaystyle{\mathop{\longrightarrow}^{I_2}} & A^{\otimes 2}&\cdots\cr\ \cr
& & \ & & \ & & \uparrow J_0^+& & \uparrow J_1^+\cr\ \cr
\ &\ & \ & \ & \ &\ & \c& 
\displaystyle{\mathop{\longrightarrow}^{id\otimes I_0}} & A\otimes\c&\cdots\cr
\ &\ & \ & \ & \ &\ & \ & 
\ & \cdots &\ 
\end{matrix}$$
where the symbols $id\otimes I_0$ and $J_0^+$ denote the unital embedding of $\c$ into $A$.
\end{lemm}

\begin{proof}
The first assertion follows from the following computation.
$$J_n^+J_n^-\begin{pmatrix}j_3 & \ldots & j_{n+1}\cr i_3 & \ldots & i_{n+1}\end{pmatrix}=
\sum \begin{pmatrix}h&g&j_3 & \ldots & j_{n+1}\cr h&g& i_3 & \ldots & i_{n+1}\end{pmatrix}=
J_n\begin{pmatrix}j_3 & \ldots & j_{n+1}\cr i_3 & \ldots & i_{n+1}\end{pmatrix}$$

The commutation of the $n$-th square in the first row follow from
$$\begin{matrix}
\sum\begin{pmatrix}h&j_2 & \ldots & j_{n}\cr h& i_2 & \ldots & i_{n}\end{pmatrix}
& \displaystyle{\mathop{\longrightarrow}^{I_{n+1}}}&
\sum\begin{pmatrix}h&j_2 & \ldots & j_{n}&l\cr h& i_2 & \ldots & i_{n}&l\end{pmatrix}
\cr &\ \cr \uparrow J_{n-1}^+& \ & \uparrow J_n^+\cr &\ \cr
\e{i_2}{j_2}\otimes \begin{pmatrix}j_3 & \ldots & j_{n}\cr i_3 & \ldots & i_{n}\end{pmatrix}
& \displaystyle{\mathop{\longrightarrow}^{id\otimes I_{n-1}}}&
\sum\e{i_2}{j_2}\otimes \begin{pmatrix}j_3 & \ldots & j_{n}&l\cr i_3 & \ldots & i_{n}&l\end{pmatrix}
\end{matrix}$$

The commutation of the $n$-th square in the second row follow from
$$\begin{matrix}
\sum\e{g}{g}\otimes\begin{pmatrix}j_3 & \ldots & j_{n+1}\cr i_3 & \ldots & i_{n+1}\end{pmatrix}
& \displaystyle{\mathop{\longrightarrow}^{id\otimes I_{n}}}&
\sum\e{g}{g}\otimes \begin{pmatrix}j_3 & \ldots & j_{n+1}&l\cr i_3 & \ldots & i_{n+1}&l\end{pmatrix}
\cr &\ \cr \uparrow J_{n}^-& \ & \uparrow J_{n+1}^-\cr &\ \cr
\begin{pmatrix}j_3 & \ldots & j_{n+1}\cr i_3 & \ldots & i_{n+1}\end{pmatrix}
& \displaystyle{\mathop{\longrightarrow}^{I_{n}}}&
\sum\begin{pmatrix}j_3 & \ldots & j_{n+1}&l\cr i_3 & \ldots & i_{n+1}&l\end{pmatrix}
\end{matrix}$$

From vertical 2-periodicity we get that the whole diagram is commutative.
\end{proof}

Define linear maps
$$E_n^-: A\otimes A^{\otimes n-1}\to A^{\otimes n-1}\hskip 20mm 
E_n^+: A^{\otimes n+1}\to A\otimes A^{\otimes n-1}$$
by the following formulae in terms of the basis.
$$E_n^-\left( \e{i_2}{j_2}\otimes \begin{pmatrix}j_3 & \ldots & j_{n+1}\cr i_3 & \ldots & i_{n+1}\end{pmatrix}\right) 
=\delta_{i_2j_2} q_{i_2}^4
\begin{pmatrix}j_3 & \ldots & j_{n+1}\cr i_3 & \ldots & i_{n+1}\end{pmatrix}$$ 
$$E_n^+\begin{pmatrix}
j_1& \ldots & j_{n+1}\cr i_1 & \ldots & i_{n+1}\end{pmatrix}
=\delta_{i_1j_1}\delta^{-2}q_{i_1}^{-4}\e{i_2}{j_2}\otimes 
\begin{pmatrix}j_3 & \ldots & j_{n+1}\cr i_3 & \ldots & i_{n+1}\end{pmatrix}$$

\begin{lemm}
The linear maps in the following diagram (E)
$$\begin{matrix}
\c &\displaystyle{\mathop{\longleftarrow}^{E_1}} & A & 
\displaystyle{\mathop{\longleftarrow}^{E_2}} & A^{\otimes 2}&
\displaystyle{\mathop{\longleftarrow}^{E_3}} & A^{\otimes 3}& 
\displaystyle{\mathop{\longleftarrow}^{E_4}} & A^{\otimes 4}&\cdots\cr\ \cr
& & \downarrow E_0^+& & \downarrow E_1^+ & & \downarrow E_2^+& & \downarrow E_3^+\cr\ \cr
\ &\ & \c & 
\displaystyle{\mathop{\longleftarrow}^{id\otimes E_0}} & A\otimes\c&
\displaystyle{\mathop{\longleftarrow}^{id\otimes E_1}} & A\otimes A& 
\displaystyle{\mathop{\longleftarrow}^{id\otimes E_2}} & A\otimes A^{\otimes 2}&\cdots\cr\ \cr
& & \ & & \downarrow E_1^- & & \downarrow E_2^-& & \downarrow E_3^-\cr\ \cr
\ &\ & \ & \ & \c&
\displaystyle{\mathop{\longleftarrow}^{I_1}} & A& 
\displaystyle{\mathop{\longleftarrow}^{I_2}} & A^{\otimes 2}&\cdots\cr\ \cr
& & \ & & \ & & \downarrow E_0^+& & \downarrow E_1^+\cr\ \cr
\ &\ & \ & \ & \ &\ & \c& 
\displaystyle{\mathop{\longleftarrow}^{id\otimes E_0}} & A\otimes\c&\cdots\cr
\ &\ & \ & \ & \ &\ & \ & 
\ & \cdots &\ 
\end{matrix}$$
are unital bimodule morphisms with respect to the inclusions in (I).
\end{lemm}

\begin{proof}
The unit for the multiplication of $A^{\otimes n}$ is
$$1_n=\sum \begin{pmatrix}l_1 & \ldots & l_n\cr l_1 & \ldots & l_n\end{pmatrix}$$

By using $(\dag)$ we get that $E_{2n}$ is unital.
$$\sum E_{2n}
\begin{pmatrix}l_1 & \ldots & l_{2n}\cr l_1 & \ldots & l_{2n}\end{pmatrix}
= \sum_{l_{2n}\sim l_{2n-1}} \delta^{-2}q_{l_{2n}}^{-4}
\begin{pmatrix}l_1 & \ldots & l_{2n-1}\cr l_1 & \ldots & l_{2n-1}\end{pmatrix}
= 1_{2n-1}$$

By using $(\ddag)$ we get that $E_{2n+1}$ is unital.
$$\sum E_{2n+1}
\begin{pmatrix}l_1 & \ldots & l_{2n+1}\cr l_1 & \ldots & l_{2n+1}\end{pmatrix}\\
= \sum q_{l_{2n+1}}^{4}
\begin{pmatrix}l_1 & \ldots & l_{2n}\cr l_1 & \ldots & l_{2n}\end{pmatrix}\\
= 1_{2n}$$

By using $(\dag )$ we get that $E_n^+$ is unital.
$$\sum E_n^+ \begin{pmatrix}l_1& \ldots & l_{n+1}\cr l_1 & \ldots & l_{n+1}\end{pmatrix}
=\sum_{l_1\sim l_2} \delta^{-2}q_{l_1}^{-4}\e{l_2}{l_2}\otimes 
\begin{pmatrix}l_3 & \ldots & l_{n+1}\cr l_3 & \ldots & l_{n+1}\end{pmatrix}=1_A\otimes 1_{n-1}$$

By using $(\ddag )$ we get that $E_n^-$ is unital.
$$\sum E_n^-\left( \e{l_2}{l_2}\otimes \begin{pmatrix}l_3 & \ldots & l_{n+1}\cr l_3 & \ldots & l_{n+1}\end{pmatrix}\right)
= \sum q_{l_2}^4\begin{pmatrix}l_3 & \ldots & l_{n+1}\cr l_3 & \ldots & l_{n+1}\end{pmatrix}=1_{n-1}$$

The right bimodule property for $E_n$ can be checked as follows.
\begin{eqnarray*}
&&E_n\left( {I_n}
\begin{pmatrix}j_1 & \ldots & j_{n-1}\cr i_1 & \ldots & i_{n-1}\end{pmatrix}
\begin{pmatrix}J_1 & \ldots & J_{n}\cr I_1 & \ldots & I_{n}\end{pmatrix}\right) \\
&=& E_n\left( \sum 
\begin{pmatrix}j_1 & \ldots & j_{n-1}&l\cr i_1 & \ldots & i_{n-1}&l\end{pmatrix}
\begin{pmatrix}J_1 & \ldots & J_{n}\cr I_1 & \ldots & I_{n}\end{pmatrix}\right) \\
&=& \delta_{(j_1\ldots j_{n-1})(I_1\ldots I_{n-1})}E_n
\begin{pmatrix}J_1 & \ldots & J_{n-1}&J_n\cr i_1 & \ldots & i_{n-1}&I_n\end{pmatrix}\\
&=& \delta_{(j_1\ldots j_{n-1})(I_1\ldots I_{n-1})}\delta_{I_nJ_n}\delta^{-1\mp 1}q_{I_n}^{\mp 4}
\begin{pmatrix}J_1 & \ldots & J_{n-1}\cr i_1 & \ldots & i_{n-1}\end{pmatrix}\\
&=& \delta_{I_nJ_n}\delta^{-1\mp 1}q_{I_n}^{\mp 4}
\begin{pmatrix}j_1 & \ldots & j_{n-1}\cr i_1 & \ldots & i_{n-1}\end{pmatrix}
\begin{pmatrix}J_1 & \ldots & J_{n-1}\cr I_1 & \ldots & I_{n-1}\end{pmatrix}\\
&=& \begin{pmatrix}j_1 & \ldots & j_{n-1}\cr i_1 & \ldots & i_{n-1}\end{pmatrix}
\wt{E}_n
\begin{pmatrix}J_1 & \ldots & J_n\cr I_1 & \ldots & I_n\end{pmatrix}
\end{eqnarray*}

The proof of the other formula $E_n(x{I_n}(y))=E_n(x)y$ is similar. For $E_n^-$ we have
\begin{eqnarray*}
&&E_n^-\left( J_n^-
\begin{pmatrix}j_3&\ldots &j_{n+1}\cr i_3&\ldots &i_{n+1}\end{pmatrix}
\left( \e{I_2}{J_2}\otimes
\begin{pmatrix}J_3&\ldots &J_{n+1}\cr I_3&\ldots &I_{n+1}\end{pmatrix}
\right)\right)\\
&=&E_n^-\left( \sum \e{g}{g}\e{I_2}{J_2}\otimes 
\begin{pmatrix}j_3&\ldots &j_{n+1}\cr i_3&\ldots &i_{n+1}\end{pmatrix}
\begin{pmatrix}J_3&\ldots &J_{n+1}\cr I_3&\ldots &I_{n+1}\end{pmatrix}\right)\\
&=& \delta_{(j_3\ldots j_{n+1})(I_3\ldots I_{n+1})}E_n^-
\left( \e{I_2}{J_2}\otimes 
\begin{pmatrix}J_3&\ldots &J_{n+1}\cr i_3&\ldots &i_{n+1}\end{pmatrix}\right)\\
&=& \delta_{(j_3\ldots j_{n+1})(I_3\ldots I_{n+1})}\delta_{I_2J_2}q_{I_2}^{4}
\begin{pmatrix}J_3&\ldots &J_{n+1}\cr i_3&\ldots &i_{n+1}\end{pmatrix}\\
&=& \delta_{I_2J_2}q_{I_2}^{4}
\begin{pmatrix}j_3&\ldots &j_{n+1}\cr i_3&\ldots &i_{n+1}\end{pmatrix}
\begin{pmatrix}J_3&\ldots &J_{n+1}\cr I_3&\ldots &I_{n+1}\end{pmatrix}\\
&=& \begin{pmatrix}j_3&\ldots &j_{n+1}\cr i_3&\ldots &i_{n+1}\end{pmatrix}
E_n^-\left( \e{I_2}{J_2}\otimes
\begin{pmatrix}J_3&\ldots &J_{n+1}\cr I_3&\ldots &I_{n+1}\end{pmatrix}
\right)
\end{eqnarray*}

The proof of the other formula $E_n^-(x{J_n^-}(y))=E_n^-(x)y$ is similar. For $E_n^+$ we have
\begin{eqnarray*}
&&E_n^+\left( J_n^+\left( \e{I_2}{J_2}\otimes
\begin{pmatrix}J_3&\ldots &J_{n+1}\cr I_3&\ldots &I_{n+1}\end{pmatrix}\right) 
\begin{pmatrix}j_1 & \ldots & j_{n+1}\cr i_1 & \ldots & i_{n+1}\end{pmatrix}\right)\\
&=& E_n^+\left( \sum 
\begin{pmatrix}h&J_2&\ldots & J_{n+1}\cr h&I_2&\ldots & I_{n+1}\end{pmatrix}
\begin{pmatrix}j_1 & \ldots & j_{n+1}\cr i_1 & \ldots & i_{n+1}\end{pmatrix}\right)\\
&=& \delta_{(J_2\ldots J_{n+1})(i_2\ldots i_{n+1})} E_n^+
\begin{pmatrix}j_1 & j_2& \ldots & j_{n+1}\cr i_1 & I_2&\ldots & I_{n+1}\end{pmatrix}\\
&=& \delta_{(J_2\ldots J_{n+1})(i_2\ldots i_{n+1})}\delta_{i_1j_1}\delta^{-2}q_{i_1}^{-4}
\e{I_2}{j_2}\otimes 
\begin{pmatrix}j_3& \ldots & j_{n+1}\cr I_3&\ldots & I_{n+1}\end{pmatrix}\\
&=& \delta_{i_1j_1}\delta^{-2}q_{i_1}^{-4}\e{I_2}{J_2}\e{i_2}{j_2}\otimes 
\begin{pmatrix}J_3&\ldots &J_{n+1}\cr I_3&\ldots &I_{n+1}\end{pmatrix}
\begin{pmatrix}j_3 & \ldots & j_{n+1}\cr i_3 & \ldots & i_{n+1}\end{pmatrix}\\
&=& \left( \e{I_2}{J_2}\otimes
\begin{pmatrix}J_3&\ldots &J_{n+1}\cr I_3&\ldots &I_{n+1}\end{pmatrix}\right)
E_n^+\begin{pmatrix}j_1 & \ldots & j_{n+1}\cr i_1 & \ldots & i_{n+1}\end{pmatrix}
\end{eqnarray*}
The proof of the other formula $E_n^+(x{J_n^+}(y))=E_n^+(x)y$ is similar.
\end{proof}

\begin{lemm}
Consider the sequence of elements $e_n$.

(i) We have $e_{2s}=J_{2s-1}J_{2s-3}\ldots J_5J_3(e_2)$ for any $s$.

(ii) We have $e_{2s+1}=J_{2s}J_{2s-2}\ldots J_6J_4J_2^+(d_2)$ for any $s$, with $d_2$ given by
$$d_2=\sum\delta^{-2}q_i^{-2}q_j^{-2}\e{i}{j}\otimes\e{i}{j}$$

(iii) The sequence $e_n$ defines a representation of the Temperley-Lieb algebra of modulus $\delta$ on the inductive limit of the algebras in the first row of (I).
\end{lemm}

\begin{proof}
By definition of $e_n$ and $J_{n-2}$ we have
$$e_n=\sum \delta^{-1\pm 1}q_i^{\pm 2}q_j^{\pm 2}\begin{pmatrix}g_1 & \ldots & g_{n-2} & j & j\cr g_1 & \ldots & g_{n-2} & i & i\end{pmatrix}=J_{n-2}(e_{n-2})$$
for any $n$. Together with
$$e_3=\sum\delta^{-2}q_i^{-2}q_j^{-2}
\begin{pmatrix}g_1&j&j\cr g_1&i&i\end{pmatrix}=J_2^+(d_2)$$
this proves (i) and (ii). By using $(\dag )$ we get that $e_{2n+1}$ is an idempotent.
\begin{eqnarray*}
e_{2n+1}^2
&=& \sum \delta^{-4}q_i^{-2}q_j^{-2}\begin{pmatrix}g_1 & \ldots & g_{2n-1} & j & j\cr g_1 & \ldots &
  g_{2n-1} & i & i\end{pmatrix}
q_I^{-2}q_J^{-2}\begin{pmatrix}G_1 & \ldots & G_{2n-1} & J & J\cr G_1 & \ldots &
  G_{2n-1} & I & I\end{pmatrix}\\
&=& \sum_{j\sim g_{2n-1}} \delta^{-4}q_j^{-4} q_i^{-2}q_J^{-2}\begin{pmatrix}g_1 & \ldots & g_{2n-1} & J & J\cr g_1 & \ldots &
  g_{2n-1} & i & i\end{pmatrix}\\
&=& e_{2n+1}
\end{eqnarray*}

By using $(\ddag )$ we get that $e_{2n+2}$ is an idempotent.

\begin{eqnarray*}
e_{2n+2}^2
&=& \sum q_i^{2}q_j^{2}\begin{pmatrix}g_1 & \ldots & g_{2n} & j & j\cr g_1 & \ldots &
  g_{2n} & i & i\end{pmatrix}
q_I^{2}q_J^{2}\begin{pmatrix}G_1 & \ldots & G_{2n} & J & J\cr G_1 & \ldots &
  G_{2n} & I & I\end{pmatrix}\\
&=& \sum q_i^{2}q_J^{2}q_j^4\begin{pmatrix}g_1 & \ldots & g_{2n} & J & J\cr g_1 & \ldots &
  g_{2n} & i & i\end{pmatrix}\\
&=& e_{2n+2}
\end{eqnarray*}

The first Jones relation can be verified as follows.
\begin{eqnarray*}
e_{3}\wt{I}_{3}(e_{2})e_{3}
&=& \sum \delta^{-4}q_i^{-2}q_k^{-2}\begin{pmatrix}K & k & k\cr K & i & i\end{pmatrix}
\sum q_h^2q_k^2\begin{pmatrix}h & h&l\cr k & k&l\end{pmatrix}
\sum q_h^{-2}q_j^{-2}\begin{pmatrix}H & j & j\cr H & h & h\end{pmatrix}\\
&=& \sum \delta^{-4}q_i^{-2}q_l^{-2}q_l^2q_l^2q_l^{-2}q_j^{-2}
\begin{pmatrix}l & j & j\cr l & i & i\end{pmatrix}\\
&=& \delta^{-2}e_{3}
\end{eqnarray*}

The other relation is proved in a similar way.
\begin{eqnarray*}
\wt{I}_{3}(e_{2})e_{3}\wt{I}_{3}(e_{2})
&=& \sum \delta^{-2}q_k^2q_i^2\begin{pmatrix}k & k&K\cr i & i&K\end{pmatrix}
\sum q_h^{-2}q_k^{-2}\begin{pmatrix}l & h & h\cr l & k & k\end{pmatrix}
\sum q_h^2q_j^2\begin{pmatrix}j & j&H\cr h & h&H\end{pmatrix}\\
&=& \sum \delta^{-2}q_l^2q_i^2q_l^{-2}q_l^{-2}q_l^2q_j^2\begin{pmatrix}j & j&l\cr i & i&l\end{pmatrix}\\
&=& \delta^{-2}\wt{I}_{3}(e_{2})
\end{eqnarray*}

By applying inclusions and shifts we get all Jones relations. 
\end{proof}

Together with $J_n^+J_n^-=J_n$ lemma 3.3 shows that the elements $e_n$ belong to the sequence of algebras obtained by going south-east starting from the algebra $A^{\otimes 2}$ in the first row of (I).

In other words, the Jones projections live at the same places as they do in standard $\lambda$-lattices axiomatized by Popa in \cite{po}.

In next four lemmas we prove that (I) together with the Jones projections and the bimodule maps in (E) satisfies Popa's axioms, namely the Jones formulae (1.1.2), the Pimsner-Popa formulae (1.3.$3^{''}$), the commuting square condition (1.1.1) and the commutation relations (2.1.1) in \cite{po}. The diagram (I) is not a standard $\lambda$-lattice is general, because the bimodule maps in (E) are not conditional expectations with respect to some trace.

\begin{lemm}
The following equalities hold
$$e_{n+2}(I_{n+2}(x))e_{n+2}=(I_{n+2}I_{n+1}E_{n+1}(x))e_{n+2}$$
$$\delta^2(I_{n+2}E_{n+2}(ye_{n+2}))e_{n+2}=ye_{n+2}$$
for any $x\in A^{\otimes n+1}$ and $y\in A^{\otimes n+2}$.
\end{lemm}

\begin{proof}
The first formula follows from the following computation.
\begin{eqnarray*}
&&e_{n+2}I_{n+2}\begin{pmatrix}j_1 & \ldots & j_{n+1}\cr i_1 & \ldots & i_{n+1}\end{pmatrix}e_{n+2}\\
&=& \sum \delta^{-1\pm 1}(q_iq_j)^{\pm 2}\begin{pmatrix}g_1 & \ldots & g_{n} & j & j\cr g_1 & \ldots &
  g_{n} & i & i\end{pmatrix}
\begin{pmatrix}j_1 & \ldots & j_{n+1}&l\cr i_1 & \ldots & i_{n+1}&l\end{pmatrix}e_{n+2}\\
&=& \sum \delta^{-1\pm 1}(q_iq_{i_{n+1}})^{\pm 2}
\begin{pmatrix}j_1 & \ldots & j_n &j_{n+1}&i_{n+1}\cr i_1 & \ldots & i_{n}&i&i\end{pmatrix}e_{n+2}\\
&=& \sum \delta^{-2\pm 2}(q_iq_{i_{n+1}}q_Iq_j)^{\pm 2}
\begin{pmatrix}j_1 & \ldots & j_n &j_{n+1}&i_{n+1}\cr i_1 & \ldots & i_{n}&i&i\end{pmatrix}
\begin{pmatrix}g_1 & \ldots & g_{n} & j & j\cr g_1 & \ldots &
  g_{n} & I & I\end{pmatrix}\\
&=& \sum \delta^{-2\pm 2}(q_iq_{i_{n+1}}^2q_j)^{\pm 2}\delta_{i_{n+1}j_{n+1}}
\begin{pmatrix}j_1 & \ldots & j_n &j&j\cr i_1 & \ldots & i_{n}&i&i\end{pmatrix}\\
&=& \sum \delta^{-2\pm 2}(q_iq_j)^{\pm 2}q_{i_{n+1}}^{\pm 4}\delta_{i_{n+1}j_{n+1}}
\begin{pmatrix}j_1 & \ldots & j_n &l_1&l_2\cr i_1 & \ldots & i_{n}&l_1&l_2\end{pmatrix}
\begin{pmatrix}g_1 & \ldots & g_{n} & j & j\cr g_1 & \ldots &
  g_{n} & i & i\end{pmatrix}\\
&=& \left( I_{n+2}I_{n+1}E_{n+1}\begin{pmatrix}j_1 & \ldots & j_{n+1}\cr i_1 & \ldots & i_{n+1}\end{pmatrix}\right) e_{n+2}
\end{eqnarray*}

The right term in the second formula is given in terms of a basis by
\begin{eqnarray*}
\begin{pmatrix}j_1 & \ldots & j_{n+2}\cr i_1 & \ldots & i_{n+2}\end{pmatrix}e_{n+2}
&=& \sum \delta^{-1\pm 1}(q_iq_j)^{\pm 2}
\begin{pmatrix}j_1 & \ldots &j_{n+2}\cr i_1& \ldots & i_{n+2}\end{pmatrix}
\begin{pmatrix}g_1 & \ldots &g_{n}&j&j\cr g_1& \ldots &g_{n}&i& i\end{pmatrix}\\
&=&\sum \delta^{-1\pm 1}(q_{j_{n+1}}q_j)^{\pm 2}\delta_{j_{n+1}j_{n+2}}
\begin{pmatrix}j_1 & \ldots &j_{n}&j&j\cr i_1& \ldots & i_n&i_{n+1}&i_{n+2}\end{pmatrix}
\end{eqnarray*}

Thus the left term is given by the following formula.
\begin{eqnarray*}
&&\delta^2I_{n+2}E_{n+2}\left( 
\begin{pmatrix}j_1 & \ldots & j_{n+2}\cr i_1 & \ldots & i_{n+2}\end{pmatrix}
e_{n+2}\right) e_{n+2}\\
&=& \sum \delta^2\delta^{-1\mp 1}\delta^{-1\pm 1}(q_{j_{n+1}}q_{i_{n+2}})^{\pm 2}\delta_{j_{n+1}j_{n+2}}q_{i_{n+2}}^{\mp 4}
\begin{pmatrix}j_1 & \ldots &j_{n}&i_{n+2}&l\cr i_1& \ldots & i_n&i_{n+1}&l\end{pmatrix}
e_{n+2}\\
&=& \sum \delta^{-1\pm 1}(q_{j_{n+1}}q_{i_{n+2}}^{-1}q_iq_j)^{\pm 2}\delta_{j_{n+1}j_{n+2}}
\begin{pmatrix}j_1 & \ldots &j_{n}&i_{n+2}&l\cr i_1& \ldots & i_n&i_{n+1}&l\end{pmatrix}
\begin{pmatrix}g_1 & \ldots &g_{n}&j&j\cr g_1& \ldots &g_{n}&i& i\end{pmatrix}\\
&=& \sum \delta^{-1\pm 1}(q_{j_{n+1}}q_{i_{n+2}}^{-1}q_{i_{n+2}}q_j)^{\pm 2}\delta_{j_{n+1}j_{n+2}}
\begin{pmatrix}j_1 & \ldots &j_{n}&j&j\cr i_1& \ldots & i_n&i_{n+1}&i_{n+2}\end{pmatrix}
\end{eqnarray*}

By cancelling $q_{i_{n+2}}^{-1}q_{i_{n+2}}$ this is equal to the right term.
\end{proof}

\begin{lemm}
The following equalities hold
$$f_{n+2}(J_{n+1}^+(x))f_{n+2}=(J_{n+1}^+J_{n+1}^-E_{n+1}^-(x))f_{n+2}$$
$$\delta^2(J_{n+1}^+E_{n+1}^+(yf_{n+2}))f_{n+2}=yf_{n+2}$$
for any $x\in A\otimes A^{\otimes n}$ and $y\in A^{\otimes n+2}$, with $f_{n+2}=I_{n+2}I_{n+1}\ldots I_4I_3(e_2)$.
\end{lemm}

\begin{proof}
The element $f_{n+2}$ is given by
$$f_{n+2}=\sum q_i^{2}q_j^{2}\begin{pmatrix}j&j&g_3 & \ldots & g_{n+2}\cr i&i&g_3 & \ldots & g_{n+2}\end{pmatrix}$$

The first formula follows from the following computation.
\begin{eqnarray*}
&&f_{n+2}J_{n+1}^+\left( \e{i_2}{j_2}\otimes
\begin{pmatrix}j_3&\ldots &j_{n+2}\cr i_3&\ldots &i_{n+2}\end{pmatrix}\right) f_{n+2}\\
&=& \sum q_i^2q_j^2
\begin{pmatrix}j&j&g_3 & \ldots & g_{n+2}\cr i&i&g_3 & \ldots & g_{n+2}\end{pmatrix}
\begin{pmatrix}h&j_2&\ldots &j_{n+2}\cr h&i_2&\ldots &i_{n+2}\end{pmatrix} f_{n+2}\\
&=& \sum q_i^2q_{i_2}^2 \begin{pmatrix}i_2&j_2&j_3&\ldots &j_{n+2}\cr i&i&i_3&\ldots &i_{n+2}\end{pmatrix} f_{n+2}\\
&=& \sum q_i^2q_{i_2}^2q_I^2q_j^2 \begin{pmatrix}i_2&j_2&j_3&\ldots &j_{n+2}\cr i&i&i_3&\ldots &i_{n+2}\end{pmatrix} 
\begin{pmatrix}j&j&g_3 & \ldots & g_{n+2}\cr I&I&g_3 & \ldots & g_{n+2}\end{pmatrix}\\
&=& \sum q_i^2q_{i_2}^4q_j^2 \delta_{i_2j_2}\begin{pmatrix}j&j&j_3&\ldots &j_{n+2}\cr i&i&i_3&\ldots &i_{n+2}\end{pmatrix}\\
&=&\sum q_i^2q_{i_2}^{4}q_j^2\delta_{i_{2}j_{2}}
\begin{pmatrix}l_1&l_2&j_3 & \ldots & j_{n+2} \cr l_1&l_2&i_3 & \ldots & i_{n+2}\end{pmatrix}
\begin{pmatrix}j&j&g_3 & \ldots & g_{n+2}\cr i&i&g_3 & \ldots & g_{n+2}\end{pmatrix}\\
&=&\left( J_{n+1}^+J_{n+1}^-E_{n+1}^-
\left( \e{i_2}{j_2}\otimes
\begin{pmatrix}j_3&\ldots &j_{n+2}\cr i_3&\ldots &i_{n+2}\end{pmatrix}\right)
\right)f_{n+2}
\end{eqnarray*}

The right term in the second formula is given in terms of a basis by
\begin{eqnarray*}
\begin{pmatrix}j_1& \ldots & j_{n+2}\cr i_1 & \ldots & i_{n+2}\end{pmatrix}f_{n+2}
&=& \sum q_i^{2}q_j^{2}
\begin{pmatrix}j_1& \ldots & j_{n+2}\cr i_1 & \ldots & i_{n+2}\end{pmatrix}
\begin{pmatrix}j&j&g_3 & \ldots & g_{n+2}\cr i&i&g_3 & \ldots & g_{n+2}\end{pmatrix}\\
&=& \sum q_{j_1}^{2}q_j^{2}\delta_{j_1j_2}
\begin{pmatrix}j&j&j_3 &\ldots & j_{n+2}\cr i_1 &i_2&i_3& \ldots & i_{n+2}\end{pmatrix}
\end{eqnarray*}

Thus the left term is given by the following formula.
\begin{eqnarray*}
&&\delta^2J_{n+1}^+E_{n+1}^+\left( 
\begin{pmatrix}j_1& \ldots & j_{n+2}\cr i_1 & \ldots & i_{n+2}\end{pmatrix}f_{n+2}\right) 
f_{n+2}\\
&=& \sum q_{j_1}^{2}q_{i_1}^{2}\delta_{j_1j_2}q_{i_1}^{-4}
\begin{pmatrix}h&i_1&j_3 &\ldots & j_{n+2}\cr h &i_2&i_3& \ldots & i_{n+2}\end{pmatrix}f_{n+2}\\
&=& \sum q_{j_1}^{2}q_{i_1}^{2}q_i^2q_j^2\delta_{j_1j_2}q_{i_1}^{-4}
\begin{pmatrix}h&i_1&j_3 &\ldots & j_{n+2}\cr h &i_2&i_3& \ldots & i_{n+2}\end{pmatrix}
\begin{pmatrix}j&j&g_3 & \ldots & g_{n+2}\cr i&i&g_3 & \ldots & g_{n+2}\end{pmatrix}\\
&=& \sum q_{j_1}^{2}q_{i_1}^{4}q_j^2\delta_{j_1j_2}q_{i_1}^{-4}
\begin{pmatrix}j&j&j_3 &\ldots & j_{n+2}\cr i_1 &i_2&i_3& \ldots & i_{n+2}\end{pmatrix}
\end{eqnarray*}

By cancelling $q_{i_{1}}^4q_{i_{1}}^{-4}$ this is equal to the right term.
\end{proof}

\begin{lemm}
The following equalities hold
$$d_{n+2}(J_{n+2}^-(x))d_{n+2}=(J_{n+2}^-J_{n}^+E_n^+(x))d_{n+2}$$
$$\delta^2(J_{n+2}^-E_{n+2}^-(yd_{n+2}))d_{n+2}=yd_{n+2}$$
for any $x\in A^{\otimes n+1}$ and $y\in A\otimes A^{\otimes n+1}$, with $d_{n+2}=(id\otimes I_{n+1}I_n\ldots I_3I_2)(d_2)$.
\end{lemm}

\begin{proof}
The element $d_{n+2}$ is given by
$$d_{n+2}=\sum \delta^{-2}q_i^{-2}q_j^{-2}\e{i}{j}\otimes
\begin{pmatrix}j&g_3 & \ldots & g_{n+2}\cr i&g_3 & \ldots & g_{n+2}\end{pmatrix}$$

The first formula follows from the following computation.
\begin{eqnarray*}
&&d_{n+2}J_{n+2}^-
\begin{pmatrix}j_2&\ldots &j_{n+2}\cr i_2&\ldots &i_{n+2}\end{pmatrix}d_{n+2}\\
&=& \sum \delta^{-2}q_i^{-2}q_j^{-2}\left(\e{g}{g} \e{i}{j}\otimes 
\begin{pmatrix}j&g_3 & \ldots & g_{n+2}\cr i&g_3 & \ldots & g_{n+2}\end{pmatrix}
\begin{pmatrix}j_2&\ldots &j_{n+2}\cr i_2&\ldots &i_{n+2}\end{pmatrix}
\right) d_{n+2}\\
&=& \sum \delta^{-2}q_{i}^{-2}q_{i_2}^{-2}\left( \e{i}{i_2}\otimes 
\begin{pmatrix}j_2&j_3&\ldots &j_{n+2}\cr i&i_3&\ldots &i_{n+2}\end{pmatrix}
\right) d_{n+2}\\
&=& \sum \delta^{-4}q_{i}^{-2}q_{i_2}^{-2} q_I^{-2}q_j^{-2}\e{i}{i_2}\e{I}{j}\otimes 
\begin{pmatrix}j_2&j_3&\ldots &j_{n+2}\cr i&i_3&\ldots &i_{n+2}\end{pmatrix}
\begin{pmatrix}j&g_3 & \ldots & g_{n+2}\cr I&g_3 & \ldots & g_{n+2}\end{pmatrix}\\
&=& \sum \delta^{-4}q_{i}^{-2}q_{i_2}^{-4}q_j^{-2} \e{i}{j}\otimes 
\begin{pmatrix}j&j_3&\ldots &j_{n+2}\cr i&i_3&\ldots &i_{n+2}\end{pmatrix}\\
&=& \sum \delta^{-4}q_{i}^{-2}q_{i_2}^{-4}q_j^{-2} \e{l_1}{l_1}\e{i}{j}\otimes
\begin{pmatrix}l_2&j_3&\ldots &j_{n+2}\cr l_2&i_3&\ldots &i_{n+2}\end{pmatrix}
\begin{pmatrix}j&g_3 & \ldots & g_{n+2}\cr i&g_3 & \ldots & g_{n+2}\end{pmatrix}\\
&=& \left( J_{n+2}^-J_{n}^+E_n^+
\begin{pmatrix}j_2&\ldots &j_{n+2}\cr i_2&\ldots &i_{n+2}\end{pmatrix}\right) d_{n+2}
\end{eqnarray*}

The right term in the second formula is given in terms of a basis by
\begin{eqnarray*}
&&\left( \e{i_1}{j_1}\otimes
\begin{pmatrix}j_2&\ldots &j_{n+2}\cr i_2&\ldots &i_{n+2}\end{pmatrix}\right) d_{n+2}\\
&=& \sum \delta^{-2}q_i^{-2}q_j^{-2}  \e{i_1}{j_1}\e{i}{j}\otimes 
\begin{pmatrix}j_2&\ldots &j_{n+2}\cr i_2&\ldots &i_{n+2}\end{pmatrix}
\begin{pmatrix}j&g_3 & \ldots & g_{n+2}\cr i&g_3 & \ldots & g_{n+2}\end{pmatrix}\\
&=& \sum \delta^{-2}q_{j_1}^{-2}q_j^{-2} \delta_{j_1j_2}\e{i_1}{j}\otimes 
\begin{pmatrix}j&j_3&\ldots &j_{n+2}\cr i_2&i_3&\ldots &i_{n+2}\end{pmatrix}
\end{eqnarray*}

Thus the left term is given by the following formula.
\begin{eqnarray*}
&&\delta^2J_{n+2}^-E_{n+2}^-\left( 
\left( \e{i_1}{j_1}\otimes
\begin{pmatrix}j_2&\ldots &j_{n+2}\cr i_2&\ldots &i_{n+2}\end{pmatrix}\right) d_{n+2}
\right) d_{n+2}\\
&=& \sum  q_{j_1}^{-2}q_{i_1}^{-2} \delta_{j_1j_2}q_{i_1}^{4}
\left( \e{g}{g}\otimes 
\begin{pmatrix}i_1&j_3&\ldots &j_{n+2}\cr i_2&i_3&\ldots &i_{n+2}\end{pmatrix}\right) d_{n+2}\\
&=& \sum \delta^{-2} q_{j_1}^{-2}q_{i_1}^{-2} \delta_{j_1j_2}q_{i_1}^{4}q_i^{-2}q_j^{-2}
\e{g}{g}\e{i}{j}\otimes 
\begin{pmatrix}i_1&j_3&\ldots &j_{n+2}\cr i_2&i_3&\ldots &i_{n+2}\end{pmatrix}
\begin{pmatrix}j&g_3 & \ldots & g_{n+2}\cr i&g_3 & \ldots & g_{n+2}\end{pmatrix}\\
&=& \sum \delta^{-2} q_{j_1}^{-2}q_{i_1}^{-4} \delta_{j_1j_2}q_{i_1}^4q_j^{-2}
\e{i_1}{j}\otimes 
\begin{pmatrix}j&j_3&\ldots &j_{n+2}\cr i_2&i_3&\ldots &i_{n+2}\end{pmatrix}
\end{eqnarray*}

By cancelling $q_{i_{1}}^{-4}q_{i_{1}}^{4}$ this is equal to the right term.
\end{proof}

\begin{lemm}
(i) The diagram obtained from (I) by replacing its vertical rows with the vertical rows of (E) commutes.

(ii) For any rectangular subdiagram of (I) having $\c$ in the south-west corner, the algebra in the north-west corner commutes with the algebra in the south-east corner.
\end{lemm}

\begin{proof}
(i) Let (IE) be this diagram. The commutation in its first row follow from
$$\begin{matrix}
\begin{pmatrix}j_1 & \ldots & j_{n}\cr i_1 & \ldots & i_{n}\end{pmatrix}
& \displaystyle{\mathop{\longrightarrow}^{I_{n+1}}}&
\sum\begin{pmatrix}j_1 & \ldots & j_{n}&l\cr i_1 & \ldots & i_{n}&l\end{pmatrix}
\cr &\ \cr \downarrow E_{n-1}^+& \ & \downarrow E_n^+\cr &\ \cr
\delta_{i_{1}j_{1}}\delta^{-2}q_{i_{1}}^{-4}\e{i_2}{j_2}\otimes \begin{pmatrix}j_3 & \ldots & j_{n}\cr i_3 & \ldots & i_{n}\end{pmatrix}
& \displaystyle{\mathop{\longrightarrow}^{id\otimes I_{n-1}}}&
\sum \delta_{i_{1}j_{1}}\delta^{-2}q_{i_{1}}^{-4}\e{i_2}{j_2}\otimes 
\begin{pmatrix}j_3 & \ldots & j_{n}&l\cr i_3 & \ldots & i_{n}&l\end{pmatrix}
\end{matrix}$$

The commutation in the second row of (IE) follow from
$$\begin{matrix}
\e{i_2}{j_2}\otimes\begin{pmatrix}j_3 & \ldots & j_{n+1}\cr i_3 & \ldots & i_{n+1}\end{pmatrix}
& \displaystyle{\mathop{\longrightarrow}^{id\otimes I_{n}}}&
\sum\e{i_2}{j_2}\otimes
\begin{pmatrix}j_3 & \ldots & j_{n+1}&l\cr i_3 & \ldots & i_{n+1}&l\end{pmatrix}
\cr &\ \cr \downarrow E_{n}^-& \ & \downarrow E_{n+1}^-\cr &\ \cr
\delta_{i_{2}j_{2}}q_{i_{2}}^{4}\begin{pmatrix}j_3 & \ldots & j_{n+1}\cr i_3 & \ldots & i_{n+1}\end{pmatrix}
& \displaystyle{\mathop{\longrightarrow}^{I_{n}}}&
\sum\delta_{i_{2}j_{2}}q_{i_{2}}^{4}\begin{pmatrix}j_3 & \ldots & j_{n+1}&l\cr i_3 & \ldots & i_{n+1}&l\end{pmatrix}\end{matrix}$$

Both diagrams (I) and (E) being 2-periodic on the vertical, (IE) is 2-periodic as well on the vertical, so it commutes.

(ii) Let $a<b$ be the ranks of the lines of (I) containing the north and south vertices of the rectangle. For $b$ odd the rectangle is of the form
$$\begin{matrix}
A^{\otimes 2s}
& \displaystyle{\mathop{\longrightarrow}^{I}}&
A^{\otimes 2s+k}& 
\cr &\ \cr \uparrow & \ & \uparrow  J\cr &\ \cr
\c
& \displaystyle{\mathop{\longrightarrow}}&
A^{\otimes k}&\end{matrix}\hskip 2cm \begin{matrix}
A\otimes A^{\otimes 2s}
& \displaystyle{\mathop{\longrightarrow}^{id\otimes I}}&
A\otimes A^{\otimes 2s+k}& 
\cr &\ \cr \uparrow & \ & \uparrow  J_{2s+k+1}^-J\cr &\ \cr
\c
& \displaystyle{\mathop{\longrightarrow}}&
A^{\otimes k}&\end{matrix}$$
depending on the parity of $a$, with $I$ and $J$ given by
$$I=I_{2s+k}I_{2s+k-1}\ldots I_{2s+2}I_{2s+1}\hskip 2cm J=J_{2s+k-1}J_{2s+k-3}\ldots J_{k+3}J_{k+1}$$

The corresponding images are
\begin{eqnarray*}
Im(I)&=&\left\{ \sum\lambda
\begin{pmatrix}j_1&\ldots &j_{2s}&l_1&\ldots &l_{k}\cr i_1&\ldots &i_{2s}&l_1&\ldots &l_{k}
\end{pmatrix}\right\}\\
Im(J)&=& \left\{ \sum\lambda
\begin{pmatrix}g_1&\ldots &g_{2s}&j_1&\ldots &j_{k}\cr g_1&\ldots &g_{2s}&i_1&\ldots &i_{k}
\end{pmatrix}\right\}
\end{eqnarray*}
and $Im(id\otimes I)=1\otimes Im(I)$ and $Im(J_{2s+k+1}^-J)=1\otimes Im(J)$, so commutation is clear. For $b$ even the rectangle is of the form
$$\begin{matrix}
A^{\otimes 2s+1}
& \displaystyle{\mathop{\longrightarrow}^{I}}&
A^{\otimes 2s+k+2}& 
\cr &\ \cr \uparrow & \ & \uparrow  J\cr &\ \cr
\c
& \displaystyle{\mathop{\longrightarrow}}&
A\otimes A^{\otimes k}&\end{matrix}\hskip 2cm 
\begin{matrix}
A\otimes A^{\otimes 2s+1}
& \displaystyle{\mathop{\longrightarrow}^{id\otimes I}}&
A\otimes A^{\otimes 2s+k+2}& 
\cr &\ \cr \uparrow & \ & \uparrow  J_{2s+k+3}^-J\cr &\ \cr
\c
& \displaystyle{\mathop{\longrightarrow}}&
A\otimes A^{\otimes k}&\end{matrix}$$
depending on the parity of $a$, with $I$ and $J$ given by
$$I=I_{2s+k+1}I_{2s+k}\ldots I_{2s+3}I_{2s+2}\hskip 2cm J=J_{2s+k+1}J_{2s+k-1}\ldots 
J_{k+5}J_{k+3}J_{k+1}^+$$

The corresponding images are
\begin{eqnarray*}
Im(I)&=&\left\{ \sum\lambda
\begin{pmatrix}j_1&\ldots &j_{2s}&j_{2s+1}&l_1&\ldots &l_{k+1}\cr i_1&\ldots &i_{2s}&i_{2s+1}&l_1&\ldots &l_{k+1}
\end{pmatrix}\right\}\\
Im(J)&=& \left\{ \sum\lambda
\begin{pmatrix}g_1&\ldots &g_{2s}&h&j_2&\ldots &j_{k+2}\cr g_1&\ldots &g_{2s}&h&i_2&\ldots &i_{k+2}
\end{pmatrix}\right\}
\end{eqnarray*}
and $Im(id\otimes I)=1\otimes Im(I)$ and $Im(J_{2s+k+3}^-J)=1\otimes Im(J)$, so commutation is clear.
\end{proof}

Define a linear form $\psi$ on $A$ by
$$\psi\e{i}{j}=\delta_{ij}p_i^{4}$$
where $p_i$ are the following positive numbers.
$$p_i=\delta^{-\frac{1}{2}}q_i^{-1}\left( \sum_{l\sim i}q_l^4\right)^\frac{1}{4}$$

By using $(\dag )$ and $(\ddag )$ we get
$$\sum p_i^4=\sum_{i\sim l}\delta^{-2}q_i^{-4}q_l^4=\sum q_l^4=1$$

This formula will be called $(\ddag )$ for $p_i$'s. There is also a corresponding $(\dag )$ formula.
$$\sum_{i\sim k}p_i^{-4}=
\sum_{i\sim k}\delta^2q_i^{4}\left( \sum_{l\sim k}q_l^4\right)^{-1}=\delta^2$$

Next lemma shows that the linear forms $\varphi_n$ define a filtered linear form $\varphi_\infty$, which fails to commute globally with vertical maps in (E) because $\varphi$ and $\psi$ are not equal in general. In fact $\varphi =\psi$ if and only if $q_i=p_i$ and the above formula for the numbers $p_i$ shows that this happens if and only if $\varphi$ is the Perron-Frobenius trace.

\begin{lemm}
Consider the diagram (I) in lemma 3.1.

(i) The linear maps $\varphi_n$ in proposition 3.1 define a filtered unital linear form $\varphi_\infty$ on the sequence of algebras in the first row of (I).

(ii) The restriction of $\varphi_\infty$ to the algebra $A\otimes A^{\otimes n-1}$ in the second row is $\psi\otimes \varphi_{n-1}$.

(iii) The diagram of restrictions of $\varphi_\infty$ is 2-periodic on the vertical.

(iv) The restrictions of $\varphi_\infty$ commute with the horizontal maps in (E).

(v) We have $\varphi_{n-1}E_n^-E_n^+=\psi_2E_{2n}$, where $E_{2n}=E_{3}E_{4}\ldots E_{n-1}E_n$.

(vi) We have $(\psi\otimes\varphi_{n-1})E_n^+=\psi_2I_2E_{1n}$, where $E_{1n}=E_{2}E_{2n}$.
\end{lemm}

\begin{proof}
(i) It is enough to verify the equality $\varphi_{n-1}E_n=\varphi_n$.
\begin{eqnarray*}
\varphi_{n-1}E_n
\begin{pmatrix}j_1 & \ldots & j_{n}\cr i_1 & \ldots & i_{n}\end{pmatrix}
&=& \delta_{i_nj_n}\delta^{-1\mp 1}q_{i_n}^{\mp 4}\varphi_{n-1}
\begin{pmatrix}j_1 & \ldots & j_{n-1}\cr i_1 & \ldots & i_{n-1}\end{pmatrix}\\
&=&  \delta_{i_nj_n}\delta^{-1\mp 1+\frac{1}{2}\pm\frac{1}{2}-n+1}q_{i_n}^{\mp 4}\delta_{(i_1\ldots i_{n-1})(j_1\ldots j_{n-1})}
q_{(i_1\ldots i_{n-1})}^4\\
&=& \delta_{(i_1\ldots i_{n})(j_1\ldots j_{n})}\delta^{\frac{1}{2}\mp\frac{1}{2}-n}
q_{(i_1\ldots i_{n})}^4\\
&=& \varphi_{n}\begin{pmatrix}j_1 & \ldots & j_{n}\cr i_1 & \ldots & i_{n}\end{pmatrix}
\end{eqnarray*}

(ii) is verified as follows.
\begin{eqnarray*}
\varphi_{n+1}J_n^+\left( 
\e{i_2}{j_2}\otimes \begin{pmatrix}j_3 & \ldots & j_{n+1}\cr i_3 & \ldots & i_{n+1}\end{pmatrix}\right)
&=& \sum\varphi_{n+1}\begin{pmatrix}h&j_2& \ldots & j_{n+1}\cr h&i_2& \ldots & i_{n+1}\end{pmatrix}\\
&=& \sum_{h\sim i_2} \delta^{\frac{1}{2}\pm\frac{1}{2}-n-1}\delta_{(i_2\ldots i_{n+1})(j_2\ldots j_{n+1})}
q_h^4q_{i_2}^{-4}q_{(i_3\ldots i_{n+1})}^4\\
&=& \delta_{i_2j_2}p_{i_2}^{4}\delta^{\frac{1}{2}\pm\frac{1}{2}-n+1}
\delta_{(i_3\ldots i_{n+1})(j_3\ldots j_{n+1})}q_{(i_3\ldots i_{n+1})}^4\\
&=& (\psi\otimes\varphi_{n-1})\left( 
\e{i_2}{j_2}\otimes \begin{pmatrix}j_3 & \ldots & j_{n+1}\cr i_3 & \ldots & i_{n+1}\end{pmatrix}\right)
\end{eqnarray*}

(iii) It is enough to show that the restrictions of $\varphi_\infty$ to the algebras in the third row are the linear forms $\varphi_n$. This follows from the equality $J_n^+J_n^-=J_n$ in lemma 3.1 and from the formulae $(\dag)$ and $(\ddag)$.
\begin{eqnarray*}
\varphi_{n+1}J_n\begin{pmatrix}j_3 & \ldots & j_{n+1}\cr i_3 & \ldots & i_{n+1}\end{pmatrix}
&=& \varphi_{n+1}
\begin{pmatrix}h&g&j_3 & \ldots & j_{n+1}\cr h&g&i_3 & \ldots & i_{n+1}\end{pmatrix}\\
&=&\sum_{g\sim h} q_h^4q_g^{-4}\delta^{-2}\varphi_{n-1}
\begin{pmatrix}j_3 & \ldots & j_{n+1}\cr i_3 & \ldots & i_{n+1}\end{pmatrix}\\
&=& \sum q_h^4\varphi_{n-1}
\begin{pmatrix}j_3 & \ldots & j_{n+1}\cr i_3 & \ldots & i_{n+1}\end{pmatrix}\\
&=& \varphi_{n-1}\begin{pmatrix}j_3 & \ldots & j_{n+1}\cr i_3 & \ldots & i_{n+1}\end{pmatrix}
\end{eqnarray*}

(iv) From proof of (i) we know that the restrictions of $\varphi_\infty$ to the algebras in the first row commute with the bimodule morphisms. By tensoring everything to the left with $id$ we get the assertion for the second row. By vertical 2-periodicity of everything this extends to the whole diagram.

(v) The map $E_{2n}$ is given by the following formula.
$$E_{2n}
\begin{pmatrix}j_1&\ldots j_n\cr i_1&\ldots i_n\end{pmatrix}
=\begin{pmatrix}j_1&j_2\cr i_1&i_2\end{pmatrix}\varphi_{n-2}
\begin{pmatrix}j_3&\ldots j_n\cr i_3&\ldots i_n\end{pmatrix}$$

The map $\varphi_{n-1}E_n^-E_n^+$ is given by 
$$\varphi_{n-1}E_n^-E_n^+
\begin{pmatrix}j_1& \ldots & j_{n+1}\cr i_1 & \ldots & i_{n+1}\end{pmatrix}
=\delta_{i_1j_1}\delta_{i_2j_2}\delta^{-2}q_{i_1}^{-4}q_{i_2}^{4}
\varphi_{n-1}\begin{pmatrix}j_3 & \ldots & j_{n+1}\cr i_3 & \ldots & i_{n+1}\end{pmatrix}$$
and together with the definition of $\psi_2$ before proposition 3.1 this proves (v).

(vi) From the defining formulae of $\psi$ and $E_n^+$ we get
\begin{eqnarray*}
(\psi\otimes\varphi_{n-1})E_n^+\begin{pmatrix}
j_1& \ldots & j_{n+1}\cr i_1 & \ldots & i_{n+1}\end{pmatrix}
&=&\delta_{i_1j_1}\delta^{-2}q_{i_1}^{-4}\psi\e{i_2}{j_2}\varphi_{n-1}
\begin{pmatrix}j_3 & \ldots & j_{n+1}\cr i_3 & \ldots & i_{n+1}\end{pmatrix}\\
&=& \delta_{i_1j_1}\delta_{i_2j_2}\delta^{-2}q_{i_1}^{-4}p_{i_2}^{4}
\varphi_{n-1}
\begin{pmatrix}j_3 & \ldots & j_{n+1}\cr i_3 & \ldots & i_{n+1}\end{pmatrix}
\end{eqnarray*}

From $i_1\sim i_2$ and from the definition of numbers $p_i$ we get
$$p_{i_1}^4q_{i_1}^4=\delta^{-2}\sum_{l\sim i_1}q_l^4=\delta^{-2}\sum_{l\sim i_2}q_l^4
=p_{i_2}^4q_{i_2}^4$$

By replacing in the above formula $q_{i_1}^{-4}p_{i_2}^{4}$ by  $q_{i_2}^{-4}p_{i_1}^{4}$ we get
$$(\psi\otimes\varphi_{n-1})E_n^+
\begin{pmatrix}j_1&\ldots &j_{n+1}\cr i_1&\ldots &i_{n+1}\end{pmatrix}
=\delta_{i_1j_1}\delta_{i_2j_2}\delta^{-2}p_{i_1}^4q_{i_2}^{-4}\varphi_{n-1}
\begin{pmatrix}j_3&\ldots &j_{n+1}\cr i_3&\ldots &i_{n+1}\end{pmatrix}$$

From the definition of numbers $p_i$ we get
$$\psi_2I_2\e{i_1}{j_1}=\sum \psi_2\ee{i_1}{j_1}{l}{l}=\sum_{l\sim {i_1}}
\delta_{i_1j_1}\delta^{-2}q_{i_1}^{-4}q_l^{4}=\delta_{i_1j_1}p_{i_1}^4$$

The map $E_{1n}$ is given by
$$E_{1n}
\begin{pmatrix}j_1&\ldots j_n\cr i_1&\ldots i_n\end{pmatrix}
=\delta_{i_2j_2}\delta^{-2}q_{i_2}^{-4}\begin{pmatrix}j_1\cr i_1\end{pmatrix}\varphi_{n-2}
\begin{pmatrix}j_3&\ldots j_n\cr i_3&\ldots i_n\end{pmatrix}$$
so the composition $\psi_2I_2E_{1n}$ is given by the same formula as $(\psi\otimes\varphi_{n-1})E_n^+$.
\end{proof}

\begin{proof}[Proof of proposition 3.1]
Let $Q_n\subset A^{\otimes n}$ be a sequence of $\c^*$-algebras satisfying conditions (1,2,3) in proposition 3.1. Define $R_n=E_{n}^+(Q_{n+1})$ and consider the following diagram ($\star$).
$$\begin{matrix}
\c &\subset & Q_{1} & \subset & Q_{2}&
\subset & Q_{3}& \subset\
\cdots\cr &\ &\ \cup &\ &\cup &\ &\cup \cr &\ &\ \c & \subset &
R_{1}&
\subset & R_{2}& \subset\ \cdots\cr &\ &\
&\ &\cup &\ &\cup \cr &\ &\ \ & \ & \c & \subset & Q_{1}&
\subset\ \cdots\cr &\ &\ &\ &\ &\ &\cup \cr
&\ &\ &\ &\ &\ &\cdots &\ \ \ \ \ \cdots\cr
\end{matrix}$$

We claim that this is a subsystem of $\c^*$-algebras of the diagram (I) in lemma 3.2.

First, the map $E_n^+$ being an involutive bimodule morphism, $R_n$ is a $\c^*$-algebra. The other thing is to verify that all inclusions in the statement make sense. The assumption $I_n(Q_{n-1})\subset Q_n$ in (1) justifies the inclusions in the first row. The bimodule property of $E_n^+$ shows that $R_n$ is included in $Q_{n+1}$ via $J_n^+$, so the first row of vertical inclusions is the good one as well. The commuting square property in lemma 3.8 (i) justifies the second row of horizontal inclusions. From condition $J_n(Q_{n-1})\subset Q_{n+1}$ in (1) and from $J_n^+J_n^-=J_n$ in lemma 3.2 and $E_n^+J_n^+=id$ in lemma 3.3 we get
$$J_n^-(Q_{n-1})=E_n^+J_n^+J_n^-(Q_{n-1})=E_n^+J_n(Q_{n-1})\subset 
E_n^-(Q_{n+1})=R_n$$
so the second row of vertical inclusions is the good one. By vertical 2-periodicity we get that ($\star$) is a subsystem of (I). 

Condition (2) and the definition of $\varphi_\infty$ show that the restrictions of $\varphi_\infty$ to the algebras in ($\star$) are traces. We prove now that ($\star$), together with $\varphi_\infty$ and with the Jones projections in lemma 3.3 is a standard $\lambda$-lattice, with $\lambda =\delta^2$.

For, it is enough to verify the commutation of $\varphi_\infty$ with the maps in (E), cf. discussion preceding lemma 3.4. Commutation with horizontal maps follows from lemma 3.8 (iv). By using vertical 2-periodicity of everything, it remains to prove that the restrictions of $\varphi_\infty$ commute with the maps of the form $E_n^+$ and $E_n^-$ in the first two rows of vertical maps of (E).

From lemma 3.8 (iv) we get $\varphi_{n+1}=\varphi_2E_{2n}=\varphi E_{1n}$. Lemma 3.8 (ii) shows that the commutation of $\varphi_\infty$ with $E_n^+$ is equivalent to the equality $(\psi\otimes\varphi_{n-1})E_n^+=\varphi_{n+1}$. This follows from lemma 3.8 (vi), from $\varphi_{n+1}=\varphi E_{1n}$ and from condition (3). The formulae for restrictions $\varphi_\infty$ in Lemma 3.8 (i,ii,iii) show that their commutation with $E_n^-$ is equivalent to the equality $\varphi_{n-1}E_n^-=\psi\otimes\varphi_{n-1}$. This must be true on the image of $E_n^+$, so is equivalent to the equality $\varphi_{n-1}E_n^-E_n^+=(\psi\otimes\varphi_{n-1})E_n^+$ obtained by composing with $E_n^+$. This follows from lemma 3.8 (v), from $\varphi_{n+1}=\varphi_2E_{2n}$ and from condition (3).

Thus ($\star$) is a standard $\lambda$-lattice with $\lambda =\delta^2$. The ``bubbling'' construction of Jones in \cite{j1} applies and proves proposition 3.1.
\end{proof}

\section{The planar algebra of a coaction -- twisted case}

Let $H$ be a Hopf $*$-algebra as in \S 1. Associated to any complex number $z$ is a multiplicative functional $f_z:H\to\c$ such that the following equalities hold (theorem 5.6 in \cite{w1}).

(f1) $f_0=\varepsilon$ and $(f_z\otimes f_t)\Delta=f_{z+t}$ for any $z,t$.

(f2) $S^2=(f_{1}\otimes id\otimes f_{-1})\Delta^{(2)}$, where $\Delta^{(2)}=(id\otimes\Delta )\Delta$.

(f3) $f_zS=f_{-z}$ and $f_z*=\wo{f}_{-\overline{z}}$ for any $z$.

(f4) $\sigma =(f_{1}\otimes id\otimes f_1)\Delta^{(2)}$ satisfies $h(ab)=h(b\sigma (a))$ for any $a,b$.

Let $A$ be a finite dimensional $\c^*$-algebra. Write $A$ as a direct sum of complex matrix algebras and let $Tr$ be the trace of $A$ which on matrix subalgebras is the usual trace of matrices.

\begin{lemm}
If $v :A\to A\otimes H$ is a coaction there exists a unique $Q\in A$ satisfying the following conditions.

(i) $ad(Q)=(id\otimes f_{\frac{1}{4}})v$.

(ii) $Q>0$.

(iii) $Tr(Q^4)=1$.

(iv) The numbers $Tr(B^{-4})$ with $B$ matrix block of $Q$ are all equal.
\end{lemm}

\begin{proof}
For any real number $z$ consider the linear map $\rho_z=(id\otimes f_z)v$. Since both $f_z$ and $v$ are multiplicative, this is an automorphism of the complex algebra $A$. By applying (f1) we get $\rho_0=id$ and
$$\rho_{z+t}=(id\otimes f_z\otimes f_t)(id\otimes\Delta)v
=(id\otimes f_z\otimes f_t)(v\otimes id)v
=(id\otimes f_z)v (id\otimes f_t)v
=\rho_{z}\rho_{t}$$

This shows that $\rho_z$ has $n$-th roots for any $n$. But $\rho_z$ must leave invariant the center $Z(A)$ of $A$, so its restriction is an automorphism of $Z(A)$ having $n$-th roots for any $n$. This is not possible if the restriction is not the identity. Since $\rho_z$ preserves the central minimal idempotents of $A$, it has to preserve the matrix blocks. But on matrix algebras automorphisms are inner, so $\rho_z$ is inner. Choose $Q_z\in A$ such that $\rho_z$ is equal to $ad(Q_z)=Q_z\, .\, Q_z^{-1}$. By using the second equality in (f3) with $z$ real we get
$$\rho_{-z}
=(id\otimes\wo{f_z})(id\otimes *)v
=*(id\otimes f_z)(*\otimes *)v
=*(id\otimes f_z)v *
=*\rho_{z}*$$
and together with $\rho_{z}\rho_{-z}=\rho_0=id$ this gives $\rho_z*\rho_{z}*=id$. On the other hand
$$\rho_z*\rho_{z}*(a)
=Q_z(Q_za^*Q_z^{-1})^*Q_z^{-1}
=ad((Q_z^*Q_z^{-1})^{-1})(a)$$
so $Q_z^*Q_z^{-1}$ is in $Z(A)$. Let $B$ be a matrix block of $Q_z$ and let $\lambda$ be a complex number such that $B^*=\lambda B$. By applying $*$ we get $B=\wo{\lambda}B^*$ and by combining these two formulae we get $B=\lambda\wo{\lambda}B$, so $\lambda$ is of modulus one. Choose a half root $\lambda^{\frac{1}{2}}$ of $\lambda$ and let $B^\prime =\lambda^{\frac{1}{2}}B$. Then
$$(B^\prime)^*=\wo{\lambda^{\frac{1}{2}}}B^*=\lambda^{-\frac{1}{2}}\lambda B=
\lambda^{\frac{1}{2}}B=B^\prime$$

Rescale in this way all blocks of $Q_z$ such that they become self-adjoint. We have $\rho_z=ad(Q_z)$ and $Q_z=Q_z^*$ for any $z$. Let $Q=Q_{\frac{1}{8}}^2$. Then $Q$ is positive and
$$ad(Q)=ad(Q_{\frac{1}{8}})^2=\rho_{\frac{1}{8}}^2=\rho_{\frac{1}{4}}=(id\otimes f_{\frac{1}{4}})v$$

We can rescale all blocks of $Q$ such that (iv) holds, then rescale $Q$ such that (iii) holds.

For the converse, if $Q^\prime$ is another element satisfying all conditions in the statement then (i) shows that $Q^\prime Q^{-1}$ is central, so if $Q=(B_i)$ is a decomposition of $Q$ then $Q^\prime$ must be of the form $(\lambda_iB_i)$. From (ii) we get that $\lambda_i >0$, then from (iv) we get that the $\lambda_i$'s are equal, and finally from (iii) we get that they are all equal to $1$.
\end{proof}

Choose a system of matrix units $X\subset A$ such that the element $Q\in A$ in lemma 4.1 is diagonal, with eigenvalues $q_i$. Let $\delta$ be the square root of the numbers in lemma 4.1 (iv). Then $\varphi =Tr(Q^4.)$ is a $\delta$-form (see \S 3).

By using boxes instead of discs as in Jones' paper \cite{j1}, we say that a tangle in $\p$ is ``vertical'' if it can be isotoped to a tangle all whose strings are parralel to the $y$-axis.

Consider the planar algebra $P(A)$ associated to the bipartite graph of the inclusion $\c\subset A$, with Perron-Frobenius spin vector. See Jones \cite{j2}.

\begin{theo}
Let $v :A\to A\otimes H$ be a coaction. Assume that $v$ preserves the linear form $\varphi =Tr(Q^4.)$ with $Q\in A$ given by lemma 4.1. There exists a unique $\c^*$-planar algebra structure $Q(v )$ on the sequence of spaces of fixed points of $v_n$ such that

(i) For any vertical tangle $T\in\p$ the multilinear map of $Q(v )$ associated to $T$ is the restriction of the multiplinear map of $P(A)$ associated to $T$.

(ii) The Jones projections are given by
$$e_n=\sum \delta^{-1\pm 1}q_i^{\pm 2}q_j^{\pm 2}\begin{pmatrix}g_1 & \ldots & g_{n-2} & j & j\cr g_1 & \ldots & g_{n-2} & i & i\end{pmatrix}$$

This $\c^*$-planar algebra is spherical and of modulus $\delta$.
\end{theo}

\begin{proof}
Uniqueness follows from the fact that $\p$ is generated by vertical tangles and Jones projections (see \S 2). It remains to verify conditions (1,2,3) in proposition 3.1.

(1) By using lemma 2.1 and lemma 2.2, it is enough to check the modularity condition in lemma 2.2. The formula of $\theta$ in \S 2 gives
$$\theta\e{i}{j}=q_i^4q_j^{-4}\e{i}{j}=ad(Q)^4\e{i}{j}=(id\otimes f_1)v\e{i}{j}$$

By using twice the axiom for coactions $(id\otimes\Delta )v =(v\otimes id)v$ we get
$$(id\otimes\Delta^{(2)})v =(id\otimes id\otimes\Delta) (v\otimes id)v 
=(v\otimes id\otimes id)(v\otimes id)v$$

By using (f4) we get that the modularity condition is satisfied.
\begin{eqnarray*}
(id\otimes\sigma )v
&=&(id\otimes f_{1}\otimes id\otimes f_{1})(v\otimes id\otimes id)(v\otimes id)v\\
&=&(((f_1\otimes id)v )\otimes id)v (id\otimes f_1)v\\
&=&(\theta\otimes id )v\theta
\end{eqnarray*}

(2) This will follow from $\theta_n =(id\otimes f_1)v_n$. In terms of the basis, the map $(id\otimes f_1)v_n$ is
$$(id\otimes f_1)v_n\begin{pmatrix}j_1 &\ldots &j_n\cr i_1 &\ldots &i_n\end{pmatrix}
=\sum \begin{pmatrix}l_1 &\ldots &l_n\cr k_1 &\ldots &k_n\end{pmatrix}
q_{(k_1\ldots k_n)}^{-1}q_{(i_1\ldots i_n)}
q_{(j_1\ldots j_n)}q_{(l_1\ldots l_n)}^{-1}$$
$$\cdot\,\, f_1V_n\begin{pmatrix}l_1&\ldots &l_n &j_1 &\ldots &j_n\cr 
k_1&\ldots &k_n &i_1 &\ldots &i_n\end{pmatrix}$$
so condition $\theta_n =(id\otimes f_1)v_n$ is equivalent to the following condition ($f_1$) for $V_n$.
$$f_1V_n\begin{pmatrix}l_1&\ldots &l_n &j_1 &\ldots &j_n\cr 
k_1&\ldots &k_n &i_1 &\ldots &i_n\end{pmatrix}
=\delta_{(_{i_1}^{j_1}\ldots {\ }_{i_n}^{j_n})(_{k_1}^{l_1}\ldots {\ }_{k_n}^{l_n})}
\,\, q_{(i_1\ldots i_n)}^4q_{(j_1\ldots j_n)}^{-4}$$

Since $\theta_1$ is the modular map of $\varphi$, this is true for $n=1$.
$$f_1V\ee{k}{l}{i}{j}=\delta_{ki}\delta_{jl}\,\, q_i^4q_j^{-4}$$

Since $f_1$ is multiplicative, from ($f_1$) for $V$ we get ($f_1$) for $V_2$.
\begin{eqnarray*}
f_1V_2\begin{pmatrix}l_1 & l_2 & j_1 & j_2 \cr k_1 & k_2 &i_1 & i_2\end{pmatrix}
&=&f_1\ee{k_1}{k_2}{i_1}{i_2}f_1\ee{l_2}{l_1}{j_2}{j_1}\\
&=&\delta_{k_1i_1}\delta_{k_2i_2}q_{i_1}^4q_{i_2}^{-4}
\delta_{l_2j_2}\delta_{l_1j_1}q_{j_2}^4q_{j_1}^{-4}\\
&=&\delta_{(k_1k_2)(i_1i_2)}\delta_{(l_1l_2)(j_1j_2)}q_{(i_1i_2)}^4q_{(j_1j_2)}^{-4}
\end{eqnarray*}

Since $f_1$ is multiplicative, from  ($f_1$) for $V$ we get ($f_1$) for $V_3$.
\begin{eqnarray*}
f_1 V_3\begin{pmatrix}l_1 & l_2 &l_3& j_1 & j_2&j_3 \cr k_1 & k_2&k_3 &i_1 & i_2&i_3\end{pmatrix}
&=&f_1 V\ee{k_1}{k_2}{i_1}{i_2}f_1V\ee{k_3}{l_3}{i_3}{j_3}f_1 V\ee{l_2}{l_1}{j_2}{j_1}\\
&=&\delta_{k_1i_1}\delta_{k_2i_2}q_{i_1}^4q_{i_2}^{-4}
\delta_{k_3i_3}\delta_{l_3j_3}q_{i_3}^4q_{j_3}^{-4}
\delta_{l_2j_2}\delta_{l_1j_1}q_{j_2}^4q_{j_1}^{-4}\\
&=&\delta_{(k_1k_2k_3)(i_1i_2i_3)}\delta_{(l_1l_2l_3)(j_1j_2j_3)}
q_{(i_1i_2i_3)}^4q_{(j_1j_2j_3)}^{-4}
\end{eqnarray*}

The proof is similar for arbitrary $n$.

(3) The coaction $v_2$ is given by the formula (see \S 1)
$$v_2\begin{pmatrix}j_1&j_2\cr i_1&i_2\end{pmatrix}
=\sum \begin{pmatrix}l_1&l_2\cr k_1&k_2\end{pmatrix}
\otimes q_{k_1}^{-1}q_{k_2}q_{i_1}q_{i_2}^{-1}q_{j_1}q_{j_2}^{-1}
q_{l_1}^{-1}q_{l_2}V\ee{k_1}{k_2}{i_1}{i_2}V\ee{l_2}{l_1}{j_2}{j_1}$$

Consider the following matrix.
$$W=\sum \begin{pmatrix}i_1&i_2\cr k_1&k_2\end{pmatrix}\otimes
q_{k_1}^{-1}q_{k_2}q_{i_1}q_{i_2}^{-1}V\ee{k_1}{k_2}{i_1}{i_2}$$

By using ($*$) we get that the matrix $W^*$ is given by
\begin{eqnarray*}
W^*
&=& \sum \begin{pmatrix}k_1&k_2\cr i_1&i_2\end{pmatrix}\otimes
q_{k_1}^{-1}q_{k_2}q_{i_1}q_{i_2}^{-1}
V\ee{k_2}{k_1}{i_2}{i_1}\\
&=& \sum \begin{pmatrix}l_1&l_2\cr j_1&j_2\end{pmatrix}\otimes
q_{l_1}^{-1}q_{l_2}q_{j_1}q_{j_2}^{-1}V\ee{l_2}{l_1}{j_2}{j_1}
\end{eqnarray*}

The linear map $x\mapsto W(x\otimes 1)W^*$ is given by
$$W\left( \begin{pmatrix}j_1&j_2\cr i_1&i_2\end{pmatrix}\otimes 1\right) W^*
=\left( \sum \begin{pmatrix}I_1&I_2\cr k_1&k_2\end{pmatrix}\otimes
q_{k_1}^{-1}q_{k_2}q_{I_1}q_{I_2}^{-1}V\ee{k_1}{k_2}{I_1}{I_2}\right)$$
$$\cdot\,\, \left( \begin{pmatrix}j_1&j_2\cr i_1&i_2\end{pmatrix}\otimes 1\right) 
\left( \sum \begin{pmatrix}l_1&l_2\cr J_1&J_2\end{pmatrix}\otimes
q_{l_1}^{-1}q_{l_2}q_{J_1}q_{J_2}^{-1}V\ee{l_2}{l_1}{J_2}{J_1}\right)$$
so we have $v_2(x)=W(x\otimes 1)W^*$ for any $x$. By using ($S$) we get
\begin{eqnarray*}
(id\otimes S)W
&=& \sum \begin{pmatrix}i_1&i_2\cr k_1&k_2\end{pmatrix}\otimes
q_{k_1}^{-1}q_{k_2}q_{i_1}q_{i_2}^{-1}SV\ee{k_1}{k_2}{i_1}{i_2}\\
&=& \sum \begin{pmatrix}i_1&i_2\cr k_1&k_2\end{pmatrix}\otimes
q_{k_1}^{-1}q_{k_2}q_{i_1}q_{i_2}^{-1}
q_{k_1}^2q_{i_1}^{-2}q_{i_2}^2q_{k_2}^{-2}
V\ee{i_2}{i_1}{k_2}{k_1}\\
&=& \sum \begin{pmatrix}i_1&i_2\cr k_1&k_2\end{pmatrix}\otimes
q_{k_1}q_{k_2}^{-1}q_{i_1}^{-1}q_{i_2}
V\ee{i_2}{i_1}{k_2}{k_1}\\
&=& \sum \begin{pmatrix}k_1&k_2\cr i_1&i_2\end{pmatrix}\otimes
q_{i_1}q_{i_2}^{-1}q_{k_1}^{-1}q_{k_2}
V\ee{k_2}{k_1}{i_2}{i_1}\\
&=& W^*
\end{eqnarray*}

By using ($\varepsilon$) and ($\Delta$) we get
$$(id\otimes\varepsilon )W=\sum \begin{pmatrix}i_1&i_2\cr k_1&k_2\end{pmatrix}
q_{k_1}^{-1}q_{k_2}q_{i_1}q_{i_2}^{-1}\delta_{k_1i_1}\delta_{k_2i_2}=\sum \begin{pmatrix}k_1&k_2\cr k_1&k_2\end{pmatrix}=1_2$$
$$(id\otimes\Delta )W=\sum \begin{pmatrix}i_1&i_2\cr k_1&k_2\end{pmatrix}\otimes 
q_{k_1}^{-1}q_{k_2}q_{i_1}q_{i_2}^{-1}V\ee{k_1}{k_2}{g}{h}\otimes V\ee{g}{h}{i_1}{i_2}$$

On the other hand $W_{12}W_{13}$ is given by
$$W_{12}W_{13}=\sum 
\begin{pmatrix}g&h\cr k_1&k_2\end{pmatrix}\begin{pmatrix}i_1&i_2\cr g&h\end{pmatrix}
\otimes q_{k_1}^{-1}q_{k_2}q_{g}q_{h}^{-1}V\ee{k_1}{k_2}{g}{h}
\otimes q_{g}^{-1}q_{h}q_{i_1}q_{i_2}^{-1}V\ee{g}{h}{i_1}{i_2}$$
and this is equal to $(id\otimes\Delta )W$. Thus $W$ is a unitary corepresentation and $v_2=ad(W)$. Consider the matrix
$$Q_W=(id\otimes f_{\frac{1}{2}})W=\sum \begin{pmatrix}i_1&i_2\cr k_1&k_2\end{pmatrix}
q_{k_1}^{-1}q_{k_2}q_{i_1}q_{i_2}^{-1}f_{\frac{1}{2}}V\ee{k_1}{k_2}{i_1}{i_2}$$

By using $(f_1)$ for $V$ we can compute $Q_W$.
$$Q_W=\sum \begin{pmatrix}i_1&i_2\cr k_1&k_2\end{pmatrix}
q_{k_1}^{-1}q_{k_2}q_{i_1}q_{i_2}^{-1}\delta_{k_1i_1}\delta_{k_2i_2}q_{i_1}^2q_{i_2}^{-2}
=\sum \begin{pmatrix}i_1&i_2\cr i_1&i_2\end{pmatrix}q_{i_1}^2q_{i_2}^{-2}$$

From the formulae (f1--f4) and from cosemisimplicity of $H$ we get that the equality $Tr(Q_W^2.)=Tr(Q_W^{-2}.)$ holds on $End(W)$ (lemma 1.1 in \cite{subf}). On the other hand, the above formula of $Q_W$ shows that $Tr(Q_W^2.)$ is $\varphi_2$ and $Tr(Q_W^{-2}.)$ is $\psi_2$. Together with the equality $End(W)=Q_2(v )$ coming from $v_2=ad(W)$, this shows that $\varphi_2=\psi_2$ on $Q_2(v )$.
\end{proof}

A natural question now is about how to verify the assumptions of theorem 4.1, namely that $v$ preserves the linear form $\varphi =Tr(Q^4.)$ with $Q\in A$ given by lemma 4.1.

(1) In the $S^2=id$ case we have $Q=1$ and the condition is satisfied. This is not very interesting, because theorem 4.1 is weaker anyway than theorem 2.1 in this case.

(2) In the $S^2\neq id$ case there are basically two examples. First is the case of adjoint coactions, where the condition is satisfied. This follows from the explicit formulae of Woronowicz in \cite{w1}, and the whole thing is discussed in detail in \cite{subf}. The other example is with the universal coactions in \cite{fc}, and once again, one can check that the condition is satisfied.

We don't know if there is a simpler characterisation of coactions $v$ for which $Q(v)$ is a planar algebra. In case there is one, getting it from what we do in this paper is probably a purely Hopf $\c^*$-algebraic problem, with no planar topology involved. This possible remaining problem is to be added to those mentioned at the end of the introduction.


\begin{thebibliography}{99}

\bibitem{subf} T. Banica, Representations of compact quantum groups and subfactors, {\em J. Reine Angew. Math.} {\bf 509} (1999), 167--198.

\bibitem{kac} T. Banica, Subfactors associated to compact Kac algebras, {\em Integral Equations Operator Theory} {\bf 39} (2001), 1--14.

\bibitem{fc} T. Banica, Quantum groups and Fuss-Catalan algebras, {\em Comm. Math. Phys.} {\bf 226} (2002), 221--232.

\bibitem{gr} T. Banica, Quantum automorphism groups of homogenous graphs, preprint math.QA/0311042.

\bibitem{bl} B. Bhattacharyya and Z. Landau, Intermediate standard invariants and intermediate planar algebras, preprint.

\bibitem{bj1} D. Bisch and V.F.R. Jones, Singly generated planar algebras
  of small dimension, {\em Duke Math. J.} {\bf 101} (2000), 41--75.

\bibitem{bj2} D. Bisch and V.F.R. Jones, Singly generated planar algebras of small dimension II, {\em Adv. Math.} {\bf 175} (2003), 297--318.

\bibitem{das} P. Das, Weak Hopf $\c^*$-algebras and depth two subfactors, preprint.

\bibitem{dak} P. Das and V. Kodiyalam, Planar algebras and the Ocneanu-Szymanski theorem, preprint.

\bibitem{da} M.-C. David, Paragroupe d'Adrian Ocneanu et alg\`ebre de Kac, {\em Pacific J. Math.} {\bf 172} (1996), 331--363.

\bibitem{ghj} F.M. Goodman, P. de la Harpe and V.F.R. Jones, ``Coxeter graphs and towers of algebras'', Springer-Verlag (1989)

\bibitem{iz} M. Izumi, Goldman's type theorems in index theory, in ``Operator algebras and quantum field theory'', Doplicher, Longo, Roberts, Zsido eds., International Press (1997), 249--271.

\bibitem{j0} V.F.R. Jones, Index for subfactors, {\em Invent. Math.} {\bf
    72} (1983), 1--25.

\bibitem{j1} V.F.R. Jones, Planar algebras, I, preprint math.QA/9909027.

\bibitem{j2} V.F.R. Jones, The planar algebra of a bipartite graph, in ``Knots in Hellas '98'', World Sci. Publishing (2000), 94--117.

\bibitem{j3} V.F.R. Jones, The annular structure of subfactors, in ``Essays on geometry and related topics'', Monogr. Enseign. Math. 38 (2001), 401--463.

\bibitem{js} V.F.R. Jones and V.S. Sunder, ``Introduction to subfactors'', {\em London. Math. Soc. Lect. Notes} {\bf 234}, Cambridge University Press (1997)

\bibitem{kls} V. Kodiyalam, Z. Landau and V.S. Sunder, The planar algebra associated to a Kac algebra, preprint.

\bibitem{la} Z. Landau, Exchange relation planar algebras, {\em Geometriae Dedicata} {\bf 95} (2002), 183--214.

\bibitem{ls} Z. Landau and V.S. Sunder, Planar depth and planar subalgebras, {\em J. Funct. Anal.} {\bf 195} (2002), 71--88.

\bibitem{lo} R. Longo, A duality for Hopf algebras and for subfactors. I, {\em Comm. Math. Phys.} {\bf 159} (1994), 133--150.

\bibitem{mvd} A. Maes and A. Van Daele, Notes on compact quantum groups, {\em Nieuw Arch. Wisk.} {\bf 16} (1998), 73--112.

\bibitem{p0} S. Popa, Classification of amenable subfactors of type II, {\em Acta Math.} {\bf 172} (1994), 163--255.

\bibitem{po} S. Popa, An axiomatization of the lattice of higher relative commutants of a subfactor, {\em Invent. Math.} {\bf 120} (1995), 427--445.

\bibitem{sa} S. Sawin, Subfactors constructed from quantum groups, {\em Amer. J. Math.} {\bf 117} (1995), 1349--1369.

\bibitem{sz} W. Szymanski, Finite index subfactors and Hopf algebra crossed products, {\em Proc. Amer. Math. Soc.} {\bf 120} (1994), 519--528.

\bibitem{tl} N.H.V. Temperley and E.H. Lieb, Relations between the ``percolation'' and ``colouring'' problem and other graph-theoretical problems associated with regular planar lattices: some exact results for the ``percolation'' problem, {\em Proc. Roy. Soc. London} {\bf 322} (1971), 251--280.

\bibitem{was} A. Wassermann, Coactions and Yang-Baxter equations for ergodic actions and subfactors, {\em London. Math. Soc. Lect. Notes} {\bf 136} (1988), 203--236.

\bibitem{w1} S.L. Woronowicz,  Compact matrix pseudogroups, {\em Comm. Math. Phys.} {\bf 111} (1987), 613--665.

\bibitem{w2} S.L. Woronowicz, Tannaka-Krein duality for compact matrix
  pseudogroups. Twisted SU(N) groups, {\em Invent. Math.} {\bf 93} (1988), 35--76.

\bibitem{w3} S.L. Woronowicz, Compact quantum groups, in ``Sym\'etries quantiques'' (Les Houches, 1995), North-Holland, Amsterdam (1998), 845--884.

\end{thebibliography}
\end{document}